\documentclass[11pt,reqno,a4paper]{article}

\usepackage{amsmath}
\usepackage{amsthm}
\usepackage{amsfonts}
\usepackage{amssymb}

\newtheorem{theorem}{Theorem}[section]
\newtheorem{lemma}[theorem]{Lemma}
\newtheorem{proposition}[theorem]{Proposition}
\theoremstyle{definition}
\newtheorem{definition}{Definition}

\newtheorem{corollary}[theorem]{Corollary}

\newcommand{\sign}{\operatorname{sign}}
\newcommand{\conv}{\operatorname{conv}}

\newcommand{\argmin}{\operatorname{arg}\min}

\newcommand{\RR}{{\mathbb{R}}}
\newcommand{\bC}{{\mathbb{C}}}
\newcommand{\bN}{{\mathbb{N}}}
\newcommand{\cK}{{\mathcal{K}}}
\newcommand{\cH}{{\mathcal{H}}}

\newtheorem{alg}{Algorithm}

\usepackage{graphics}
\usepackage{graphicx}
\begin{document}
\title{Recovery algorithms for vector valued data with joint sparsity constraints}
\author{Massimo Fornasier\footnotemark[2] and Holger Rauhut\footnotemark[3]}
\date{August 2006}
\maketitle
\begin{abstract}
Vector valued data appearing in concrete applications often possess sparse 
expansions with respect to a preassigned frame for each vector component individually. 
Additionally, different components may also exhibit common sparsity patterns. 
Recently, there were introduced 
sparsity measures that take into account such joint sparsity patterns, promoting coupling of 
non-vanishing components. These measures are typically constructed as weighted $\ell_1$ norms of 
componentwise $\ell_q$ norms of frame coefficients. We show how to compute solutions of linear inverse 
problems with such joint sparsity regularization constraints by fast thresholded Landweber algorithms. 
Next we discuss the adaptive choice of suitable weights appearing in the definition of sparsity measures. 
The weights are interpreted as indicators of the sparsity pattern and are iteratively up-dated after each 
new application of the thresholded Landweber algorithm. The resulting two-step algorithm is interpreted 
as a double-minimization scheme for a suitable target functional. We show its $\ell_2$-norm convergence.
An implementable version of the algorithm is also formulated, and its norm convergence is proven. 
Numerical experiments in color image restoration are presented.
\end{abstract}



\renewcommand{\thefootnote}{\fnsymbol{footnote}}

\footnotetext[2]{Johann Radon Institute for Computational and Applied Mathematics,
Austrian Academy of Sciences,\\ 
Altenbergerstrasse 69, A-4040 Linz, Austria, 
{\tt massimo.fornasier@oeaw.ac.at}\\
MF acknowledges the financial support provided by the
European Union's Human Potential Programme under contract MEIF-CT-2004-501018. He also thanks
NuHAG for its warm hospitality.}
\footnotetext[3]{Numerical Harmonic Analysis Group, 
Faculty of Mathematics, University of Vienna,\\
Nordbergstrasse 15, A-1090 Vienna, Austria,
{\tt holger.rauhut@univie.ac.at}\\ 
HR acknowledges the financial support provided by the
European Union's Human Potential Programme under contracts HPRN--CT--2002--00285 (HASSIP)
and MEIF-CT-2006-022811.}

\renewcommand{\thefootnote}{\arabic{footnote}}

\section{Introduction}

{\bf Inverse problems}. 
We address the problem of recovering 
an element $u$ of a Hilbert space $\mathcal{K}$ from the observed 
datum $g = T u$ in the Hilbert space $\mathcal{H}$, where 
$T: \mathcal{K} \rightarrow \mathcal{H}$ is a bounded linear operator, possibly
non-invertible or with unbounded inverse.
A simple approach to this problem is
to minimize the discrepancy
$$
        \mathcal{T}(u):= \| T u - g|{\mathcal{H}} \|^2.
$$
If $\text{ker}(T)=\{0\}$ then there exists a unique solution given by $u^* = (T^* T)^{-1} T^* g$. However,
if $T$ has unbounded inverse, i.e., $(T^*T)^{-1}$ is unbounded then this approach is very unstable.

Thus, if $T$ is non-invertible or has unbounded inverse (or an inverse with high norm)
one has to take into account further features of the expected solution.
Indeed, a well-known way out is to consider the regularized problem \cite{enhane96}
\begin{eqnarray*}
        u^*_\alpha&:=& \text{argmin}_{u \in \mathcal{K}} \mathcal{T}(u) + \alpha \|u|\mathcal{K}\|^2.
\end{eqnarray*} 
for which the corresponding solution operator $T_\alpha^\dagger: g \mapsto u^*_\alpha$ is bounded. 
Unfortunately, the minimal norm constraint is often not appropriate.
A recent approach is to substitute this particular constraint with a more general one 
\begin{eqnarray*}
        u^*_\Phi&:=& \text{argmin}_{u \in \mathcal{K}} \mathcal{T}(u) + \Phi(u),
\end{eqnarray*} 
where $\Phi$ is a suitable {\it sparsity measure}. 

{\bf Sparse frame expansions}. A {\it sparse representation} of an element of a Hilbert space
is a series expansion with respect to an orthonormal
bases or a frame that has only a small number of large coefficients.
%
Several types of signals appearing in nature admit 
sparse frame expansions and thus, sparsity is a 
realistic assumption for a very large class of problems. For instance, images are well-represented
by sparse expansions with respect to wavelets or curvelets, while for audio signals 
a Gabor frame is a good choice.

Sparsity has had already a long history of successes. The design of frames 
for sparse representations of digital signals  has led to extremely efficient compression methods, 
such as JPEG2000 and MP3 
\cite{ma99}. Successively a new generation of optimal numerical schemes has been developed for the 
computation of sparse solutions of differential and integral equations, exploiting adaptive and 
greedy strategies \cite{coh03,S,DFR,DFRSW}. The use of sparsity in inverse problems 
for data recovery has been the most recent step of this long career of 
``simplifying and understanding complexity'', with an enormous 
potential in applications \cite{an05,dadede04, date05, do92,do95,do95-1,rate05,te05,caur04, cohore04,dama03}.
Another field, which caught much attention recently, is the observation
that it is possible to reconstruct sparse signals from vastly incomplete information
\cite{cataXX,carotaXX,do04,kura06,ra05-7}.
This line of research is called sparse recovery or compressed sensing.

{\bf From sparsity to joint sparsity}. Most of the contributions appearing in the literature 
are addressed to the recovery of sparse scalar functions.
Multi-channel signals (i.e., vector valued functions) appearing in concrete applications 
may not only possess sparse frame expansions for each channel individually, 
but additionally the different channels can also exhibit common sparsity patterns.
Recently, new sparsity measures have been introduced that promote such coupling
of the non-vanishing components through different channels \cite{babadusawa05,gisttr06,tr06}.
These measures are typically constructed as {\it weighted} $\ell_1$ norms of channel $\ell_q$ norms
with $q>1$.
We will use this concept for the solution of vector valued inverse problems and
combine it with another approach further promoting the coupling of sparsity patterns
along channels.

{\bf Our main results}. We show how to compute solutions of linear 
inverse problems with joint sparsity 
regularization constraints by fast thresholded Landweber algorithms, similar to those presented in \cite{dadede04,rate05, te05}. We discuss the adaptive choice of 
suitable {\it weights} appearing in the definition of the sparsity measures. The {\it weights} are 
interpreted as indicators of the sparsity pattern and are iteratively up-dated after 
each new application of the thresholded Landweber algorithm. The resulting two-step algorithm is 
interpreted as a double-minimization scheme for a suitable target functional. 
We prove that our algorithm converges to its minimizer.
Since the functional is not smooth, this is done by subdifferential inclusions \cite{rowe98}.
We prove that the thresholded Landweber algorithm, which constitutes the inner iteration of the double-minimization 
algorithm, converges linearly. This feature was not ensured by the versions in \cite{dadede04,te05}. 
The second step of the double-minimization has actually a simple explicit solution.
Finally, we show that the full exact double-minimization scheme converges linearly and we 
provide an implementable version which is also ensured to converge. 

{\bf Morphological analysis of signals and sparsity patterns}. 
The use of sparseness measures not only allows to reconstruct
a signal. At the same time it gives information about (joint)
sparsity patterns which may encode morphological features of the signal.
Well-known examples are the microlocal analysis properties of 
wavelets \cite{ja04,ja05} for singularity and regularity detection, 
and the characterization of 
edges and curves by curvelets for natural images \cite{cado04}. 
For instance, the weight sequences appearing in the sparsity measures we define, 
and interpreted as indicators of the sparsity pattern, play a similar role as the discontinuity 
set is playing in the Mumford-Shah functional \cite{mush89}. In fact, it is well-known that wavelet
or curvelet coefficients have high absolute values at high 
scales as soon as we are in the neighborhood 
of discontinuities. 
Even more illuminating  and suggestive is the parallel between the 
sparsity measure and its 
indicator weights with the Ambrosio-Tortorelli \cite{amto90} 
approximation of the Mumford-Shah functional.
Here, the discontinuity set is indicated by an auxiliary function which is 1 
where the image is smooth and 0 where edges and discontinuities are detected. 

Joint sparsity patterns of vector 
valued (i.e., multi-channel) signals encode even finer 
properties of the morphology which do not belong only to one channel but are a common feature of all 
the channels. Here the parallel is with generalizations of the Mumford-Shah functional 
as appearing for example in \cite{brkiso03} where polyconvex functions of gradients 
couple discontinuity sets through different color channels of images.

{\bf Applications}. We expect that our scheme can be applied in 
several different contexts. In this paper we limit ourselves to 
an application in color image reconstruction, modeling a real-world problem in art restoration.
Indeed, color images have the advantage to be non-trivial multivariate and multi-channel signals, 
exhibiting a very rich morphology and structure. 
In particular, discontinuities (jump sets) may appear in all the channels at the same locations, which
will be reflected in their curvelet representation (for instance).
For these reasons, color images are a good model to test the 
effectiveness of our scheme promoting joint 
sparsity, also because the solution can be easily checked just by a visual analysis. 
Of course, the range of applicability of our approach is not limited to color image restoration. 
Neuroimaging (functional Magnetic Resonance Imaging, Magnetoencephalography), distributed compressed sensing
\cite{babadusawa05} 
and several other problems with coupled vector valued solutions are fields 
where we expect that our scheme can be fruitfully used.
The numerical solution of differential and integral operator equations can also be addressed 
within this framework and we refer for example to \cite{DFR,S,DFRSW} for implementations by adaptive strategies.

{\bf Content of the paper}. The paper is organized as follows. In Section 2 we introduce the 
mathematical setting. We formulate our model of joint sparsity for multi-channel signals and 
the corresponding functional to be minimized in order to solve a given linear inverse problem. 
The functional depends on two variables. The first belongs to the space of signals 
to be reconstructed, the second belongs to the space of sparsity indicator weights. 
Convexity properties of the functional are discussed. 
Section 3 is dedicated to the formulation of the double-minimization algorithm and to 
its weak-convergence. The scheme is based on alternating minimizations in the first and 
in the second variable individually. In Section 4 we discuss an efficient thresholded 
Landweber algorithm for the minimization with respect to the first variable.  Its strong convergence 
is shown following the analysis in \cite{dadede04}. The minimization with respect to 
the second variable has an explicit solution and no elaboration is needed. We provide 
an implementable version of the full scheme in Section 5. To prove its convergence we develop an 
error analysis. As a byproduct of the results in this section we show that the 
double--minimization scheme converges strongly.
In Section 6 we present an application in color image reconstruction. 
Numerical experiments are shown and discussed.

\subsection*{Nota on color pictures}

This paper introduces methods to recover colors in digital images. Therefore a gray level 
printout of the manuscript does not allow to appreciate fully the quality of the illustrated 
techniques. The authors recommend the interested reader to access the electronic version with 
color pictures which is available online.

\section{The Functional}

\subsection{Notation}

Before starting our discussion let us briefly introduce some of the spaces
we will use in the following. For some countable index set $\Lambda$
we denote by $\ell_p=\ell_p(\Lambda)$, $1 \leq p \leq \infty$, 
the space of real sequences $u=(u_\lambda)_{\lambda \in \Lambda}$ 
with norm
\[
\|u\|_p \,=\, \|u|\ell_p\| \,:=\, 
\left(\sum_{\lambda \in \Lambda} |u_\lambda|^p\right)^{1/p}, \quad 1\leq p < \infty
\]
and $\|u\|_\infty \,:=\, \sup_{\lambda \in \Lambda} |u_\lambda|$ as usual.
If $(v_\lambda)$ is a sequence of positive weights then we define the weighted
spaces $\ell_{p,v} = \ell_{p,v}(\Lambda) = \{u, (u_\lambda v_\lambda) \in \ell_p(\Lambda)\}$
with norm
\[
\|u\|_{p,v} \,=\, \|u|\ell_{p,v}\| \,=\, \|(u_\lambda v_\lambda)\|_p \,=\,
\left(\sum_{\lambda \in \Lambda} v_\lambda^p |u_\lambda|^p)\right)^{1/p}
\]
(with obvious modification for $p=\infty$). If the entries $u_\lambda$ are actually
vectors in a Banach space $X$ with norm $\|\cdot\|_X$ then we denote
\[
\ell_{p,v}(\Lambda,X) \,:=\, \{(u_\lambda)_{\lambda \in \Lambda}, u_\lambda \in X, 
(\|u_\lambda\|_X)_{\lambda \in \Lambda} \in \ell_{p,v}(\Lambda)\}
\]
with norm $\|u|\ell_{p,v}(\Lambda,X)\| = \|(\|u_\lambda\|_X)_{\lambda \in \Lambda}|\ell_{p,v}(\Lambda)\|$.
Usually $X$ will be $\RR^M$ endowed with the Euclidean norm, or the $M$-dimensional
space $\ell_q^M$, i.e., $\RR^M$ endowed with the $\ell_q$-norm. By $\RR_+$ we denote the non-negative
real numbers.

\subsection{Inverse Problems with joint sparsity constraints}

Let $\cK$ and $\cH_j$, $j=1,\hdots,N$, be (separable) Hilbert spaces 
and $A_{\ell,j}: \cK \to \cH_j$, $j=1,\hdots,M$, $\ell=1,\hdots,N$, some bounded 
linear operators. 
Assume we are given data $g_j \in \cH_j$,
\[
g_j \,=\, \sum_{\ell = 1}^M A_{\ell,j} f_\ell,\quad j=1,\hdots,N.
\]
Then our basic task consists
in reconstructing the (unknown) elements $f_\ell \in \cK$, $\ell=1,\hdots,M$.

In practice, it happens that the corresponding mapping from the vector
$(f_\ell)$ to the vector $(g_j)$ is not invertible or ill-conditioned. Moreover,
the data $g_j$, $j=1,\hdots,N$, are often corrupted by noise. 
Thus, in order to deal with our reconstruction problem we have to regularize it.

Our basic assumption throughout this paper will be that the 'channels' 
$f_\ell$, $\ell=1,\hdots,M$, 
are correlated by means of joint sparsity patterns.
Our aim is to model the joint sparsity within a regularization term. 
In the following we develop this idea.

For the sake of short notation we resume the data vector into
\[
g \,=\, (g_j)_{j=1,\hdots,M} \in \cH:= \bigoplus_{j=1}^N \cH_j
\]
where the Hilbert space $\cH$ is equipped with the usual inner product 
$\langle \sum_j g_j, \sum_j h_j\rangle := \sum_j \langle g_j, h_j\rangle$ with $g_j,h_j \in \cH_j$.
We also combine the operators $A_{\ell,j}$ into one operator
\[
A: \bigoplus_{\ell=1}^M \cK \to \cH,\quad A (f_\ell)_{\ell=1}^{M} 
\,=\, \left(\sum_{\ell=1}^M A_{\ell,j} f_\ell\right)_{j=1}^N.
\] 

In order to exploit sparsity ideas we assume that we have given a suitable frame 
$\{\psi_\lambda: \lambda \in \Lambda\} \subset \cK$ indexed by a countable set $\Lambda$. 
This means
that there exist constants $C_1,C_2 > 0$ such that
\begin{equation}\label{frame_ineq}
C_1 \|f\|^2_\cK \leq \sum_{\lambda \in \Lambda} |\langle f, \psi_\lambda\rangle |^2 \leq C_2 
\|f\|_\cK^2\qquad \mbox{ for all } f \in \cK.
\end{equation}
Orthonormal bases are particular examples of frames. Frames allow for a (stable) series expansion
of any $f \in \cK$ of the form
\begin{equation}\label{expand_f}
f \,=\, F u \,:=\, \sum_{\lambda \in \Lambda} u_\lambda \psi_\lambda
\end{equation}
where $u = (u_\lambda)_{\lambda \in \Lambda} \in \ell_2(\Lambda)$. 
The linear operator $F : \ell_2(\Lambda) \to \cK$ is called the {\it synthesis map} in frame theory.
It is bounded due to the frame inequality (\ref{frame_ineq}).
In contrast to orthonormal bases, the coefficients $u_\lambda$ need not be unique, in general.
For more information on frames we refer to \cite{ch03-1}.

A main assumption here is that the $f_\ell$ to be reconstructed are sparse 
with respect to the frame $\{\psi_\lambda\}$.
This means that $f_\ell$ can be well-approximated by a series of the form (\ref{expand_f}) with 
only a small number of non-vanishing coefficients $u_\lambda$. This can be modelled by assuming
that the sequence $u$ is contained in a (weighted) $\ell_1(\Lambda)$-space. 
Indeed, the minimization of the $\ell_1(\Lambda)$ norm promotes that only few entries are non-zero.
Taking for instance a wavelet frame and a suitable weight, the $\ell_1$ constraint implies that 
the element to be reconstructed
lies in a certain Besov space $B_{1,1}^s$, see \cite{dadede04}.

Analogously as in
\cite{dadede04} such considerations lead to the regularized functional
\begin{align}\label{def_orig_func}
\mathcal{J}(u) \,=\, \|g - Tu|\cH\|^2 + \|u|\ell_{1,v}(\Lambda^M)\|
\,=\, \sum_{j=1}^N \left\|g_j - \sum_{\ell=1}^M A_{\ell,j} F u^{\ell}\right\|^2_{\cH_j} + 
\sum_{\ell=1}^M \sum_{\lambda \in \Lambda} v_\lambda |u_\lambda^{\ell}|,
\end{align}
which has to be minimized with respect to the vector of coefficients 
$u= (u_\lambda^{\ell})_{\lambda \in \Lambda}^{\ell = 1,\hdots,M}$. 
The $\ell_{1,v}$ norm in this functional clearly represents the regularization term. 
The numbers $v_\lambda$, $\lambda \in \Lambda$,
are some suitable positive weights. Once the minimizer $u = (u_\lambda^{\ell})$ is determined 
we obtain a reconstruction of
the vectors of interest by means of $f_\ell \,=\, F u^{\ell} = \sum_\lambda u^{\ell}_\lambda \psi_\lambda$. 
The algorithm in \cite{dadede04}
can be taken to perform the minimization with respect to $u$. 

The functional $\mathcal{J}(u)$ in the form stated, however, does not necessarily model 
any correlation 
between the vectors ('channels') $f_\ell$, $\ell=1,\hdots,M$. 
A way to incorporate such correlation
is the assumption of joint sparsity, see also \cite{gisttr06,tr06}. 
By this we mean 
that the pattern of non-zero coefficients
representing $f_\ell$ is (approximately) the same for all the channels. 
In other words, for some
{\it finite} set of indexes $\Lambda_0\subset \Lambda$ and for all $\ell=1,\hdots,N$ there 
is an expansion
\[
f_\ell\, \approx\, \sum_{\lambda \in \Lambda_0} u^{\ell}_\lambda \psi_\lambda.
\]
In particular, the {\it same} $\Lambda_0$ can be chosen for all $f_\ell$'s.


We propose two approaches (that can be combined) 
to model 
joint sparsity. 
The first one assumes that the mixed
norm 
\[
 \|u|\ell_{1,v}(\Lambda, \ell_q^M)\|
\,=\, \sum_{\lambda \in \Lambda} v_\lambda \|u_\lambda\|_q
\]
of $u=(u_\lambda^{\ell})$ is small. 
Hereby, $u_\lambda$ denotes the vector 
$(u_\lambda^{(\ell)})_{\ell=1}^M$ in $\RR^M$. (Recall also that $\ell_q^M$ denotes
$\RR^M$ endowed with the $\ell_q$-norm). 
Here, $q > 1$ and in particular, $q=2$ or $q= \infty$, represent
the interesting cases, since for $q=1$ the above norm reduces to the
usual weighted $\ell_{1,v}$ norm. 
In fact if $q$ is large and some $|u_\lambda^\ell|$ is large then the channel 
entries $|u_\lambda^{\ell'}|$ are also allowed to be large for $\ell' \neq \ell$, 
without increasing significantly the norm $\|u_\lambda|\ell_q^M\|$. 
The minimization of the above norm promotes that all entries of the 'interchannel' 
vector $u_\lambda$ may become significant, once at least one of the components 
$|u_\lambda^\ell|$ is large.

Introduce the operators $T_{\ell,j} = A_{\ell,j} F : \ell_2(\Lambda) \to \cH_\ell$ and 
\[
T : \ell_2(\Lambda,\RR^M) \to \cH,\quad T u \,=\, \left(\sum_{\ell=1}^M T_{\ell,j} u^{\ell} \right)_{j=1}^N
\,=\,  \left(\sum_{\ell=1}^M A_{\ell,j} F u^{\ell} \right)_{j=1}^N.   
\]
The above reasoning leads to the functional
\begin{align}\label{def_Psi0}
K(u) \,=\, K_{v}^{(q)}(u) \,&:=\, \|T u - g|\cH\|^2 + \|u|\ell_{1,v}(\Lambda,\ell_q^M)\|\\
&=\, \sum_{j=1}^N \left\|\sum_{\ell=1}^M T_{\ell,j} u^{\ell} - g_j\right\|^2_{\cH_j}
+ \sum_{\lambda \in \Lambda} v_\lambda \|u_\lambda\|_q \notag
\end{align} 
to be minimized with respect to $u$. 
In Section \ref{sec_iter} we will develop an iterative thresholding 
algorithm similar as in \cite{dadede04} to perform
this minimization.

The second approach to support joint sparsity is to encode the joint 
sparsity information in some sort of indicator function. This can in fact be
done by using the weight $(v_\lambda)$ as a second
minimization variable. To this end we add an additional term
to the original functional (\ref{def_orig_func}), punishing
small values of $v_\lambda$. We obtain the functional
\[
J_0(u,v) \,:=\, J_{\theta,\rho,0}^{(1)}(u,v)\,:=\,
\|T u - g|\cH\|^2 + \sum_{\lambda \in \Lambda} v_\lambda \|u_\lambda\|_1 
+ \sum_{\lambda} \theta_\lambda(\rho_\lambda - v_\lambda)^2 
\]
restricted to $v_\lambda \geq 0$. Here, $(\theta_\lambda)_\lambda$ 
and $(\rho_\lambda)_\lambda$ are some
suitable positive sequences.

Now the task is to minimize $J_0(u,v)$ 
jointly with respect to both $u,v$. 
(Again, once this minimizer is determined we obtain 
$f_\ell = F u^{\ell}$). Analyzing $J_0(u,v)$ we realize
that for the minimizer $(u,v)$ we will have $v_\lambda = 0$ (or close to $0$)
if $\|u_\lambda\|_1 = \sum_{\ell=1}^M |u_\lambda^{\ell}|$ 
is large so that $v_\lambda \|u_\lambda\|_1$ gets small. On the other hand,
if $\|u_\lambda\|_1$ is small then the term 
$\theta_\lambda(\rho_\lambda - v_\lambda)$ dominates and 
forces $v_\lambda$ to be close to $\rho_\lambda$. 
Thus, $v_\lambda$ serves indeed as an indicator of large values
of $\|u_\lambda\|_1$. It has the effect, that if $v_\lambda$ is chosen small 
due to one large $u_\lambda^\ell$ then also the other coefficients
$u_\lambda^{\ell'}$, $\ell' \neq \ell$ can be chosen large
without making the functional considerably bigger.



Unfortunately, in contrast to the previous functionals, 
$J(u,v)$ as stated above is no longer jointly convex in $(u,v)$ in general
(although it is convex as functional of $u$ and of $v$ alone). 
Thus, it cannot be ensured that a local minimum of the functional
will be a global one, a property that is very crucial for an efficient
minimization method.

To overcome this problem we may add an additional suitable
quadratic term. Moreover, we can, of course, combine the second
approach with the first one and use an $\ell_q$-norm instead of
an $\ell_1$-norm for the 'interchannel'
vectors $u_\lambda$. This leads to the most general form of the 
regularized functional considered in this paper,
\begin{align}\label{def_func_J}
J(u,v) \,=\, J^{(q)}_{\theta,\rho,\omega}(u,v) 
\,:=\, \|Tu -g|\cH\|^2 + \sum_{\lambda \in \Lambda} v_\lambda \|u_\lambda\|_q
+ \sum_{\lambda \in \Lambda} \omega_\lambda \|u_\lambda\|_2^2
+ \sum_{\lambda \in \Lambda} \theta_\lambda(\rho_\lambda - v_\lambda)^2.
\end{align}
Here, $\omega_\lambda$ is a suitably chosen sequence of positive numbers, and
$1\leq q \leq \infty$.

We will provide a sufficient condition depending on
$\theta_\lambda$ and $\rho_\lambda$ in the next subsection
ensuring the strict joint convexity of $J(u,v)$ in $(u,v)$. 
Although there is an extra term, $J(u,v)$ has similar properties 
as $J_0(u,v)$. In particular, $v$ can still be seen as a sort
of indicator function.

Observe that in the minimum we will always have $0 \leq v_\lambda \leq \rho_\lambda$.
Therefore, we can assume the domain of $J$ to be 
$\ell_2(\Lambda,\RR^M) \times \ell_{\infty,\rho^{-1}}(\Lambda)_+$
where $\ell_{\infty,\rho^{-1}}(\Lambda)_+$ denotes the (convex) 
cone of all non-negative 
sequences $(v_\lambda) \in \ell_{\infty,\rho^{-1}}(\Lambda)$. 

Our main contribution consists in providing an algorithm for
the minimization of $J(u,v)$. It consists in alternately minimizing
with respect to $u$ and with respect $v$. The minimization with respect
to $v$ can be done explicitly. 
For the minimization with respect to $u$
we propose an efficient iterative algorithm.


We will mainly study the problem in the real-valued case. The complex-valued case
can be treated with the same methods (in principle) by observing
that $\bC^M$ is isomorphic to $\RR^{2M}$, so passing from
$M$ complex-valued channels to $2M$ real-valued channels. We note, however, that slight
complications may arise from the fact that an $\ell_q$ norm on $\bC^M$ is not isometric
to an $\ell_q$-norm on $\RR^{2M}$ if $q\neq 2$. (In particular, the thresholding
operator on $\bC^M$ for $q=\infty$ will have a different form than
the one provided in the next Section for the real-valued case).

\subsection{Convexity of the functional $J$}

At several places in the following it will be convenient to write
\begin{equation}\label{def_TPhi}
J(u,v) \,=\, \mathcal{T}(u) + \Phi^{(q)}(u,v)
\end{equation}
where
\begin{align}
\mathcal{T}(u) \,&=\, \|Tu - g\|_\cH^2 \,=\, 
\sum_{j=1}^N \|\sum_{\ell=1}^M T_{\ell,j} u^{\ell} - g_j\|_{\mathcal{H}_j}^2 \notag\\
\Phi^{(q)}(u,v) 
\, &=\, 
\sum_{\lambda \in \Lambda} v_\lambda \| u_\lambda\|_q +  \sum_{\lambda \in \Lambda} \omega_\lambda \| u_\lambda\|_2^2 + \sum_{\lambda\in \Lambda} \theta_\lambda(\rho_\lambda - v_\lambda)^2 \notag\\
\,&=\, \left \| (v_\lambda \| u_\lambda\|_q)_{\lambda  \in \Lambda} \right \|_1 +  
\left \| u| \ell_{2,\omega^{1/2}}(\Lambda,\ell_2^M)\right\|^2 + \|\rho-v|\ell_{2,\theta^{1/2}}(\Lambda)\|^2,\notag
\end{align}
are the discrepancy with respect to the data and the joint sparsity measure, respectively.

Also it is useful to observe that $\Phi^{(q)}$ decouples with respect to $\lambda$, i.e.,
\begin{align}\label{sumphi}
\Phi^{(q)}(u,v) \,=\, \sum_{\lambda \in \Lambda} \Phi_\lambda^{(q)}(u_\lambda,v_\lambda)
\end{align}
where
\begin{align}\label{def_Philambda}
\Phi^{(q)}_\lambda(x,y) \,&=\, y \|x\|_q + \omega_\lambda \|x\|_2^2 + 
\theta_\lambda(\rho_\lambda - y)^2 \\
&=\, y \left(\sum_{\ell = 1}^M |x_\ell|^q\right)^{1/q}
+ \omega_\lambda \sum_{\ell=1}^M x_\ell^2 + \theta_\lambda(\rho_\lambda - y)^2, 
\qquad x \in \RR^M, y \geq 0 \notag
\end{align}
(with the usual modification for $q= \infty$).

In the following we give necessary and sufficient conditions
for the (strict) convexity of the functional $\Phi^{(q)}$ for
the most interesting cases $q=1,2,\infty$. These imply sufficient
conditions for the (strict) convexity of $J=J^{(q)}_{\theta,\rho,\omega}$.



\begin{proposition}
\label{q12oo} Let $q \in \{1,2,\infty\}$.
The sparsity measure $\Phi^{(q)}$ is convex if and only if  
$\omega_\lambda \theta_\lambda \geq \frac{\kappa}{4}$ for 
all $\lambda \in \Lambda$, where $\kappa= M$ for $q=1$, and $\kappa =1$ for 
$q \in \{2,\infty\}$. In particular, if $\omega_\lambda \theta_\lambda \geq \frac{\kappa}{4}$ 
for all $\lambda \in \Lambda$ then $J$ is convex. In case of a strict
inequality $\omega_\lambda \theta_\lambda > \frac{\kappa}{4}$ we can replace ``convexity''
by ``strict convexity'' in all of these statements.
\end{proposition}
\begin{proof}
It is easy to see that $\Phi^{(q)}$ is (strictly) convex if and only if all the
$\Phi^{(q)}_\lambda$, $\lambda \in \Lambda$, are (strictly) convex.

Let us first consider $q=1$. 
Observe that we can write $\Phi^{(1)}_\lambda(x,y) = \sum_{\ell=1}^M F^{(1)}_{\lambda}(x_\ell,y)$
with
\begin{align}\label{def_F1}
F^{(1)}_{\lambda}(z,y) \,&=\, 
y |z| + \omega_\lambda |z|^2 + M^{-1} \theta_\lambda (\rho_\lambda-y)^2
\\
\,&=\, \left( y|z| + \omega_\lambda |z|^2 + M^{-1} \theta_\lambda y^2 \right) + 
\left(\rho_\lambda ^2 -  2\theta_\lambda y \right),
\quad z \in \RR, y \geq 0.\notag
\end{align}
The function in the second bracket is obviously linear, hence convex. The function in the first bracket
can be written as the composition $G_{\lambda} \circ L$ with $L(z,y) = (|z|,y)$ and 
\[
G_{\lambda}(z,y) \,=\, yz + \omega_\lambda z^2 + M^{-1}\theta_\lambda y^2
\,=\, \frac{1}{2} (z,y) \,H^{(1)}\, (z,y)^T,\qquad z,y \in \RR, 
\]
where
\[
H^{(1)} \,=\, \left(\begin{matrix}
2 \omega_\lambda & 1 \\
1& 2 \theta_\lambda M^{-1}
\end{matrix} \right ).
\]
Since $L$ is convex and has range $\RR_+^2$, and
$G_\lambda$ is monotonically increasing in each coordinate on $\RR_+^2$ it suffices
to show that $G_{\lambda}(z,y)$ is convex if and only if $\theta_\lambda \omega_\lambda \geq M/4$,
see e.g.~\cite[p.~86]{bova04}. The convexity of $G_\lambda$ is equivalent
to $H^{(1)}$ being positive semidefinite. The latter is clearly equivalent
to $\theta_\lambda \omega_\lambda \geq M/4$. Strict convexity is equivalent to
a strict inequality $\theta_\lambda \omega_\lambda > M/4$.

Now let $q=2$. Observe that $\Phi^{(2)}_\lambda = F^{(2)}_\lambda \circ L^{(2)}$
where $L^{(2)} : \RR^M \times \RR_+ \to \RR_+^2$, $L^{(2)}(x,y) \,=\, (\|x\|_2, y)$
and
\begin{align}
F^{(2)}_\lambda(z,y) = y z + \omega_\lambda z^2 + \theta_\lambda (\rho_\lambda-y)^2
\,=\, \frac{1}{2}(z,y)\,H^{(2)}\,(z,y)^T + \theta_\lambda(\rho_\lambda^2 - 2\rho_\lambda y),
\quad z,y \in \RR
\end{align}
with
\[
H^{(2)} \,=\, \left(\begin{matrix}
2 \omega_\lambda & 1 \\
1& 2 \theta_\lambda
\end{matrix} \right ).
\]
By a similar argument as above $\Phi^{(2)}_\lambda$ is convex if and only
if $H^{(2)}$ is positive semi-definite. The latter is the case if and only
if $\omega_\lambda \theta_\lambda \geq 1/4$, and strict convexity is equivalent
to a strict inequality.


Finally, let $q=\infty$.
%
Observe that
\[
        \Phi^{(\infty)}_\lambda(x,y)  = 
\max_{\ell=1,...,M} \left \{ y |x_\ell| + \omega_\lambda \sum_{m=1}^M |x_m|^2 + \theta_\lambda (\rho_\lambda - y)^2 \right \}.
\]
Since $\Phi^{(\infty)}_\lambda$ is the pointwise 
maximum of $M$ functions, it is 
sufficient (see \cite[p.~80]{bova04}) to investigate the (strict) convexity of each of the functions
\begin{align}
f_{\lambda,\ell}(x,y) \,&=\, y x_\ell + \omega_\lambda \sum_{m=1}^M (x_\ell)^2 
+ \theta_\lambda (\rho_\lambda - y)^2\notag\\
\,&=\, \frac{1}{2} (y,x)\, H^{(\infty)}_\ell \, (y,x)^T 
+ \theta_\lambda(\rho_\lambda^2 - 2y\rho_\lambda),\quad x\in \RR^M,y\geq 0, 
\notag
\end{align}
with
$$
        H^{(\infty)}_\ell \,=\, \left(\begin{matrix}
2 \omega_\lambda & 0& \cdots &0& \delta_{1,\ell} \\
0 & 2 \omega_\lambda & \cdots &0& \delta_{2,\ell} \\
\vdots & \vdots & \ddots & \vdots & \vdots \\
0&0& \cdots &2 \omega_\lambda&\delta_{M,\ell}\\
\delta_{1,\ell} & \delta_{2,\ell}& \cdots &\delta_{M,\ell}&2 \theta_\lambda
\end{matrix} \right ).
$$
One can show by induction that
$$
\det(H^{(\infty)}_\ell) \,=\, 2^{M-1} \omega_\lambda^{M-1}
(4 \theta_\lambda \omega_\lambda -1).
$$
Thus, $H^{(\infty)}_\ell$ is positive semidefinite if and only if
$\theta_\lambda \omega_\lambda \geq 1/4$, and the convexity of 
$\Phi^{(q)}_\lambda$ is equivalent to the latter condition. Once again
strict convexity is equivalent to the strict inequality.
\end{proof}





We do not pursue the task to obtain conditions for
the convexity of $\Phi^{(q)}$ and $J$ for general $q \neq 1,2,\infty$, 
but rather assume that $\Phi^{(q)}$ and hence $J$ are always convex
also in this case.

\section{The Minimizing Algorithm and its Convergence}

In this section we propose and analyze an algorithm for the 
computation of the minimizer $(u^*,v^*)$ of the functional 
$J(u,v) = J^{(q)}_{\theta,\rho,\omega}(u,v)$ defined in (\ref{def_func_J}).
The algorithm consists in alternating a minimization with respect to $u$
and a minimization with respect to $v$.
More formally, for some initial choice $v^{(0)}$, 
for example $v^{(0)} = (\rho_\lambda)_{\lambda \in \Lambda}$,
we define
\begin{equation}
\label{firstalg}
\begin{array}{ll}
u^{(n)} := \text{arg} \min_{u \in \ell_{2}(\Lambda,\RR^M)} J(u,v^{(n-1)}),\\
v^{(n)} := \text{arg} \min_{v \in \ell_{\infty,\rho^{-1}}(\Lambda)_+} J(u^{(n)},v).
\end{array}
\end{equation}
The minimization of $J(u,v^{(n-1)})$ with respect to $u$ can be done
by means of the iterative thresholding algorithm that we will study in the next section. The minimizer $v^{(n)}$ 
of $J(u^{(n)},v)$ for fixed $u^{(n)}$ can be computed explicitly. 
Indeed, it follows from elementary calculus that
\begin{equation}\label{min_v}
v_\lambda^{(n)} \,=\, \left \{ \begin{array}{ll} 
\rho_\lambda - \frac{1}{2 \theta_\lambda} \|{u^{(n)}}_\lambda\|_q
& \mbox{ if } \|{u^{(n)}}_\lambda\|_q < 2 \theta_\lambda \rho_\lambda\\
0& \text{ otherwise }.
\end{array}
\right.
\end{equation}

We have the following result about the convergence of the above algorithm.

\begin{theorem}\label{thm_conv} 
Let $1\leq q \leq \infty$ and assume that $\Phi^{(q)}$ and
hence $J$ are strictly convex (see also Proposition \ref{q12oo}).
Moreover, we assume that $\ell_{2,\omega^{1/2}}(\Lambda,\RR^M)$ is
embedded into $\ell_{2}(\Lambda,\RR^M)$, i.e., $\omega_\lambda \geq \gamma >0$ for all $\lambda \in \Lambda$.
Then the sequence $(u^{(n)},v^{(n)})_{n \in \mathbb{N}}$ converges to 
the unique minimizer 
$(u^*,v^*) \in \ell_2(\Lambda,\RR^M) \times \ell_{\infty,\rho^{-1}}(\Lambda)_+$ 
of $J$. The convergence of $u^{(n)}$ is weak in $\ell_2(\Lambda,\RR^M)$ and that of $v^{(n)}$ holds componentwise.\\
For the most interesting cases $q \in \{1,2,\infty\}$, if in addition 
$\theta_\lambda \omega_\lambda \geq \sigma > \phi_q/4$  for all $\lambda \in \Lambda$, 
where $\phi_1 = M$, $\phi_2 = 1$, $\phi_\infty = \sqrt{M}$
then the convergence of $u^{(n)}$ to $u^*$ is also strong in  $\ell_{2}(\Lambda,\mathbb{R}^M)$ 
and $v^{(n)}-v^*$ converges to $0$ strongly in $\ell_{2,\theta}(\Lambda)$.
\end{theorem}

The rest of the section will be spent with the proof of the weak convergence of the algorithm. 
The strong convergence and the full proof of the Theorem \ref{thm_conv} will be established 
only in Subsection 5.3 later.

\subsection{Subdifferential calculus}

A main tool in the analysis of non-smooth functionals and their
minima is the concept of subdifferential. Recall that 
for a convex functional $F$ 
on some Banach space $V$ its subdifferential $\partial F(x)$
at a point $x \in V$ with $F(x) < \infty$ is defined as the set
\[
\partial F(x) \,=\, \{x^* \in V^*, x^*(z-x) + F(x) \,\leq\, F(z) 
\mbox{ for all } z \in V\},
\]
where $V^*$ denotes the dual space of $V$. It is obvious from this
definition that $0 \in \partial F(x)$ if and only if $x$ is a minimizer
of $F$. In the following we investigate the subdifferential of $J$.
In order to have $J$ defined on the whole Banach space
$\ell_2(\Lambda,\RR^M) \times \ell_{\infty,\rho^{-1}}(\Lambda)$ rather than just for 
positive $v_\lambda$'s (which is needed
to use subdifferentials) we simply 
extend $J(u,v)$ by 
\[
J(u,v) \,=\, \infty \qquad \mbox{ if } v_\lambda < 0 \mbox{ for some } \lambda
\in \Lambda
\]
as usual. This extension preserves convexity
and does not change the minimizer.

Recall that $J$ can be written as $J(u,v) = \mathcal{T}(u) + \Phi^{(q)}(u,v)$,
see (\ref{def_TPhi}). Since both $\mathcal{T}$ and $\Phi^{(q)}$ are convex
we have, see e.g.~\cite[Proposition 5.6]{ekte99},
\begin{equation}\label{subdiff_J}
\partial J(u,v) \,=\, \partial \mathcal{T}(u) \times \{0\} + 
\partial \Phi^{(q)}(u,v).
\end{equation}
Concerning the subdifferential of $\mathcal{T}$ we have the following result.

\begin{lemma}\label{lem_subdiffT}
The subdifferential of $\mathcal{T}$ at $u \in \ell_2(\Lambda,\RR^M)$ consists of one element,  
$$
\partial \mathcal{T}(u) = \left \{ 2 T^* (T u -g) \right \}.
$$
\end{lemma}
\begin{proof}
Since $\mathcal{T}$ is convex and Gateaux-differentiable, by Proposition 5.3 \cite{ekte99} we have
$
\partial \mathcal{T}(u) = \left \{ \mathcal{T}'(u) \right \},   
$
where its Gateaux-derivative is characterized by 
$\langle \mathcal{T}'(u),z \rangle = \lim_{h \rightarrow 0^+} \frac{\mathcal{T}(u+h z)-\mathcal{T}(u)}{h}$ 
for all $z \in \ell_2(\Lambda,\RR^M)$. 
It is straightforward to check that the Gateaux derivative of a functional
of the type $u \rightarrow \| T u - g\|^2$ (with linear $T$) at $u$ applied on $z$ is given by
$
2 \langle Tu - g, Tz \rangle \,=\, 2\langle T^*(Tu - g),z \rangle.
$
This proves the claim.
\end{proof}

Let us now consider the subdifferential of $\partial \Phi^{(q)}(u,v)$. Recall its
domain $\ell_2(\Lambda,\RR^M) \times \ell_{\infty,\rho^{-1}}(\Lambda)$.
Since the dual of $\ell_{\infty,\rho^{-1}}$ is a bit inconvenient to handle
we restrict the subdifferential to the predual $\ell_{1,\rho}$. This will be
enough for our purposes.
Moreover, recall that $\Phi^{(q)}$ decouples
into a sum of functionals $\Phi^{(q)}_\lambda$ depending only $(u_\lambda, v_\lambda)$, see
(\ref{sumphi}). 
It is straightforward to show the following lemma.

\begin{lemma} The subdifferential of $\Phi^{(q)}$ at a point 
$(u,v) \in \ell_2(\Lambda,\RR^M) \times \ell_{\infty,\rho^{-1}}(\Lambda)$ 
with $\Phi^{(q)}(u,v) < \infty$ 
satisfies
\begin{align}
D \Phi^{(q)}(u,v) \,&:=\, \partial \Phi^{(q)}(u,v) \cap (\ell_2(\Lambda,\RR^M) \times \ell_{1,\rho}(\Lambda))\notag\\
\,&=\, 
\{ (\xi,\eta) \in \ell_2(\Lambda,\RR^M) \times \ell_{1,\rho}(\Lambda):~ 
(\xi_\lambda,\eta_\lambda) \in \partial \Phi^{(q)}_\lambda(u_\lambda,v_\lambda)\mbox{\rm { }for all }\lambda\in \Lambda \}.
\notag
\end{align}
\end{lemma}

We are left with investigating the subdifferential of the functional 
$\Phi^{(q)}_\lambda$ defined in (\ref{def_Philambda}). Similarly as $J$ we 
extend it to $\RR^M \times \RR$ by 
$\Phi^{(q)}_\lambda(x,y) = \infty$ for  $y < 0$.
\begin{lemma} 
Let $1\leq q \leq \infty$. Assume that $\Phi^{(q)}_\lambda$ is convex 
(see also Proposition \ref{q12oo}). Then for 
$(x,y) \in \RR^M \times \RR_+$ we have
\begin{align}
\partial \Phi^{(q)}_\lambda(x,y) \,=\, 
\{(\xi,\eta) \in \RR^M \times \RR:~ \xi \in y \partial \|\cdot\|_q(x) +2 \omega_\lambda x,~ 
\label{subdiff_Phi}
\eta \in  \|x\|_q \partial s^+(y)  + 2 \theta_\lambda (y-\rho_\lambda)\}.
\end{align}
where $s^+(y) := y$ for $y \geq 0$ and $s^+(y) = \infty$ for $y < 0$. In particular,
$\partial s^+(y) = \{1\}$ for $y > 0$ and $\partial s^+(0) = (-\infty,1]$.
\end{lemma}

\begin{rem} We recall that the subdifferential of the $q$-norm on $\RR^M$ 
is given as follows.
If $1 < q < \infty$ then
\[
\partial \|\cdot\|_q(x) \,=\, 
\left\{ \begin{array}{ll} B^{q'}(1)& \mbox{ if } x = 0,\\
\left\{
\left(\frac{|x_\ell|^{q-1} \sign(x_\ell)}{\|x\|_q^{1-1/q}}\right)_{\ell=1}^M \right\} & \mbox{ otherwise},
\end{array}
\right. 
\]
where $B^{q'}(1)$ denotes the ball of radius $1$ in the dual norm, i.e., in $\ell_{q'}$ with 
$1/q+1/{q'} = 1$. 

If $q=1$ then 
\begin{equation}\label{subdiff_l1}
\partial \|\cdot\|_1(x) \,=\, \{ \xi \in \RR^M:~ \xi_\ell \in \partial |\cdot|(x_\ell), \ell=1,\hdots,M\}
\end{equation}
where $\partial |\cdot|(z) = \{\sign(z)\}$ if $z \neq 0$ and $\partial |\cdot|(0) = [-1,1]$.

If $q= \infty$ then 
\begin{equation}\label{max_subgrad}
\partial \|\cdot\|_\infty(x) \,=\, \left\{\begin{array}{ll}
B^1(1) & \mbox{if } x = 0,\\
\conv\{(\sign(x_\ell) e_\ell: |x_\ell| = \|x\|_\infty\} & \mbox{otherwise},
\end{array}\right.  
\end{equation}
where $\conv A$ denotes the convex hull of a set $A$ and $e_\ell$ the $\ell$-th canonical
unit vector in $\RR^M$. 
\end{rem}

\begin{proof} Recall that 
\[
\Phi^{(q)}_\lambda(x,y) \,=\,  s^+(y) \|x\|_q + \omega_\lambda \|x\|_2^2 + 
\theta_\lambda(\rho_\lambda - y)^2.
\]
Let $y \geq 0$ so that $\Phi^{(q)}_\lambda(x,y)$ is finite.
The subdifferential $\partial (\Phi^{(q)}_\lambda)_x(x,y)$ of $\Phi^{(q)}(x,y)$ considered as
a function of $x$ alone (i.e. for fixed $y$) is clearly given by
\begin{equation}\label{subdiff_Phix}
\partial (\Phi^{(q)}_\lambda)_x(x,y) \,=\, y \partial \|\cdot\|_q(x) + 2\omega_\lambda x 
\end{equation}
while keeping $y$ fixed gives
\[
\partial (\Phi^{(q)}_\lambda)_y(x,y) \,=\, \partial s^+(y) \|x\|_q + 2\theta_\lambda(y-\rho_\lambda).
\]
This shows the inclusion '$\subset$' in (\ref{subdiff_Phi}). Moreover, for all the points 
$(x,y) \in \RR^M \times \RR_+$ where $\Phi^{(q)}_\lambda$ is differentiable we even have equality
in (\ref{subdiff_Phi}) since $\Phi^{(q)}_\lambda$ is convex and, thus, 
all the subdifferentials appearing consist of precisely one point, i.e., the usual gradient. 

Let $1 < q < \infty$. Then for
$x \neq 0$, $y > 0$ the 
differentiability assumption is clearly satisfied. 
For the other cases $x = 0$ or $y = 0$ we note that by convexity
of $\Phi^{(q)}_\lambda$ we have (see \cite[Corollary 10.11]{rowe98})
\begin{equation}\label{subdiff_part}
\partial (\Phi_\lambda^{(q)})_x(x,y) \,=\, \{\xi: \exists \eta \mbox{ such that } (\xi,\eta) \in 
\partial \Phi^{(q)}_\lambda(x,y) \}
\end{equation}
and the corresponding relation for $\partial (\Phi_\lambda^{(q)})_y(x,y)$.
Now, if $y > 0$ then $\Phi_\lambda^{(q)}(x,y)$ is differentiable with respect
to $y$ and thus, $\eta$ in the right hand side of (\ref{subdiff_part}) is 
unique, indeed $\eta = \eta_0 := \frac{\partial}{\partial y} \Phi^{(q)}_\lambda(x,y)$.
We conclude that for $y>0$
\[
\partial(\Phi^{(q)}_\lambda)(x,y) = \{(\xi,\eta_0), \xi \in \partial 
(\Phi^{(q)}_\lambda)_x(x,y)\}
\]
In particular this holds for $x = 0$, even
for general $1\leq q \leq \infty$.
The same argument applies for the case $y = 0$ 
and $x \neq 0$ (and $1 < q < \infty$), which shows (\ref{subdiff_Phi}) in these cases.
Now let $x=0$ and $y=0$. Then the right hand side of (\ref{subdiff_Phi})
contains precisely one point, i.e., $(\xi,\eta)=(0,-2\theta_\lambda \rho_\lambda)$.
Since the subdifferential $\Phi_\lambda^{(q)}(0,0)$ contains at least one point
by convexity, it must coincide with $(\xi,\eta)$ by the trivial
inclusion '$\subset$'. (It is easy to check also directly 
that $(0,-2\theta_\lambda \rho_\lambda) \in \Phi_\lambda^{(q)}(0,0)$).
%
Note that this argument applies also for $q=1,\infty$.

It remains to treat the cases $q=1,\infty$ with $x \neq 0$ and arbitrary $y\geq 0$.
Let us start with $q=1$. In the proof of Proposition \ref{q12oo} it was noted that
\[
\Phi^{(1)}_\lambda(x,y) \,=\, \sum_{\ell=1}^M F_\lambda^{(1)}(x_\ell,y)
\]
with $F_\lambda^{(1)} : \RR^2 \to \RR$ defined in (\ref{def_F1}). The subdifferential of $F_\lambda$ 
can be obtained in the same way as above (expressing e.g. formally the modulus
as a $2$-norm on $\RR^1$). For $(z,y) \in \RR \times \RR_+$ this yields
\[
\partial F_\lambda^{(1)}(z,y) \,=\, \{(\tau,\eta):~ \tau \in y \partial |\cdot|(z) + 2\omega_\lambda z,
\eta \in |z| \partial s^+(y) + 2M^{-1} \theta_\lambda(y-\rho_\lambda) \}.
\]
By convexity we have
\[
\partial \Phi^{(1)}_\lambda(x,y) \,=\, \sum_{\ell=1}^M \left\{ (e_\ell z_\ell,\eta):~ (z_\ell,\eta) \in 
\partial F_\lambda^{(1)}(x_\ell,y)\right\}
\]
where $e_\ell$ denotes the $\ell$-th unit vector in $\RR^M$.
By the explicit form of the subdifferential of the $\ell_1$-norm
(\ref{subdiff_l1}) this gives (\ref{subdiff_Phi}) for $q=1$.

Finally, let $q=\infty$. Similarly as in the proof of Proposition \ref{q12oo} we write
\[
\Phi^{(\infty)}_\lambda(x,y) \,=\, \max_{\ell=1,\hdots,M} F_{\ell}(x,y)
\] 
with
\[
F_\ell(x,y) \,=\, y|x_\ell| + \omega_\lambda \|x\|_2^2 + \theta_\lambda(\rho_\lambda - y)^2.
\]
If $x_\ell \neq 0$ then $F_\ell(x,y)$ is differentiable with respect to $x$ and
\[
\partial F_\ell(x,y) \, = \, \{ (\xi,\eta): \xi \,=\, 
y \sign(x_\ell)e_\ell + 2\omega_\lambda x,\,y \in \partial s^+(y) |x_\ell| +2\theta_\lambda (y-\rho_\lambda)
\},
\]
where $e_\ell$ denotes the $\ell$-th canonical unit vector in $\RR^M$.
This even holds for $y = 0$ by an analogous argument as above, see (\ref{subdiff_part}).
The subdifferential of $\Phi^{(\infty)}_\lambda(x,y)$ for $x\neq 0$ is then given by (see e.g.
\cite[Exercise 8.31]{rowe98})
\[
\partial \Phi^{(\infty)}_\lambda(x,y) \,=\, \conv \{ \partial F_\ell(x,y): 
F_\ell(x,y) \,=\, \max_{m=1,\hdots,M} F_m(x,y)\}.
\]
Since $x \neq 0$ we have $x_\ell \neq 0$ if $|x_\ell| = \|x\|_\infty$ and 
the latter is the case iff $F_\ell(x,y) = \max_m F_m(x,y)$. Thus, we obtain
\begin{align}
&\partial \Phi^{(\infty)}_\lambda(x,y)\notag\\ 
&=\, \conv \bigcup_{\ell: |x_\ell| = \|x\|_\infty} 
\{ (\xi,\eta):~ \xi = 
y \sign(x_\ell)e_\ell + 2\omega_\lambda x,\eta \in \partial s^+(y) |x_\ell| +2\theta_\lambda (y-\rho_\lambda)
\}\notag\\
&=\, 
\left\{ (\xi,\eta):~ \xi \in \conv\{ y \sign(x_\ell) e_\ell,|x_\ell| = \|x\|_\infty\},
 \eta \in \|x\|_\infty \partial s^+(y) \} 
+ (2\omega_\lambda x, 2\theta_\lambda(y-\rho_\lambda)\right\} .\notag
\end{align}
By the characterization of the subdifferential of the $\infty$-norm in 
(\ref{max_subgrad}) we obtain the claimed equality in (\ref{subdiff_Phi}) for $q=\infty$ 
and $x \neq 0$. This finishes the proof. 
\end{proof}

Combining the previous lemmas 
we obtain the following result.
\begin{proposition}
\label{subdiff} Let $1\leq q \leq \infty$.
Assume that $\Phi^{(q)}$ is convex and let 
$(u,v) \in \ell_2(\Lambda,\RR^M) \times \ell_{\infty,\rho^{-1}}(\Lambda)$ such that
$\Phi^{(q)}(u,v) < \infty$. Then we have
\begin{align}
D \Phi^{(q)}(u,v) \,&= \,
\{(\xi,\eta) \in \ell_2(\Lambda,\RR^M) \times \ell_{1,\rho}(\Lambda),
\xi_\lambda \in v_\lambda \partial\|\cdot\|_q(u_\lambda) + 
2\omega_\lambda u_\lambda, \notag\\
& \phantom{:=\, \{(\xi,\eta) \in \ell_2(\Lambda,\RR^M) \times \ell_{1,\rho}(\Lambda),}
\eta_\lambda \in \|u_\lambda\|_q \partial s^+(v_\lambda) + 2\theta_\lambda(v_\lambda - \rho_\lambda),
~\lambda \in \Lambda \} \notag\\
\label{def_D}
&\subset\, \partial \Phi^{(q)}(u,v)
\end{align}
and
\[
DJ(u,v) \,=\, \partial J(u,v) \cap \left(\ell_2(\Lambda,\RR^M) \times \ell_{1,\rho}(\Lambda)\right)
\,=\, (2 T^*T(u-g),0) + D \Phi_\lambda^{(q)}(u,v) \,\subset\, \partial J(u,v).
\]
\end{proposition}

\subsection{Weak convergence of the double-minimization}

Before we actually start proving the weak convergence of 
the algorithm in (\ref{firstalg}) we recall the following definition 
\cite{rowe98}.
\begin{definition}
Let $V$ be a topological space and $\mathcal{A}=(A_n)_{n \in \mathbb{N}}$ 
a sequence of subsets of $V$. The subset $A \subseteq V$ is 
called the \emph{limit of the sequence $\mathcal{A}$}, and 
we write $A = \lim_n A_n$, if 
$$
        A=\{a \in V: \exists a_n \in A_n, a = \lim_{n} a_n\}. 
$$
\end{definition}
The following observation will be useful for us, see e.g.~\cite[Proposition 8.7]{rowe98}. 
\begin{lemma}\label{lem_osc}
Assume that $\Gamma$ is a convex function on 
$\RR^M$ and $(x_n) \subset \RR^M$ a convergent sequence 
with limit $x$ such that $\Gamma(x_n),\Gamma(x) < \infty$.
Then the subdifferentials satisfy   
$$
        \lim_n \partial \Gamma(x_n) \subseteq \partial \Gamma(x).
$$
In other words, the subdifferential $\partial \Gamma$ of a convex function
is an outer semicontinuous set-valued function. 
\end{lemma}
In the following we agree on the convention that the upper index $n$ at 
$u^{(n)} \in \ell_2(\Lambda,\RR^M)$ always denotes the $n$-th iterate and
$u_\lambda^{(n)} \in \RR^M$ denotes the (vector-valued) entry at $\lambda$
of the $n$-th iterate. In the following proof we will never refer to
the $\ell$-th component of the $M$-dimensional vector $u_\lambda^{(n)}$, so hopefully 
no confusion can arise. Also, we denote by 
$(DJ(u,v))_\lambda$ 
the restriction of $DJ(u,v) \subset \ell_2(\Lambda,\RR^M) \times \ell_{1,\rho}(\Lambda)$ 
to the index $\lambda$.
By the previous section it holds
\begin{equation}\label{DJ_lambda}
(DJ(u,v))_\lambda \,=\, (2T^*T(u-g))_\lambda,0) + \partial \Phi^{(q)}_\lambda(u_\lambda,v_\lambda).
\end{equation}   


Now the proof is developed as follows.
First, we recall that $(u^*,v^*) = \argmin J(u,v)$ if and 
only if $0 \in \partial J(u^*,v^*)$. Next, we show that 
there exist weakly convergent subsequences of $(u^{(n)}, v^{(n)})$ 
(again denoted by  $(u^{(n)}, v^{(n)})$) which converge 
to $(u^{(\infty)},v^{(\infty)})$ and that 
\begin{equation}
\label{firstdiffinc}
        0 \in \lim_n D J(u^{(n)},v^{(n)}) \subseteq \partial J(u^{(\infty)},v^{(\infty)}).
\end{equation}
Due to the strict convexity of $J$ we conclude that 
$(u^{(\infty)},v^{(\infty)}) = (u^*,v^*)$. Now, let us detail the argument.

By definition of $u^{(n)}$ and $v^{(n)}$ we have
\begin{align}
& J(u^{(n)},v^{(n)})- J(u^{(n+1)},v^{(n+1)})\notag\\ 
=\, & J(u^{(n)},v^{(n)})- J(u^{(n+1)},v^{(n)})
+ J(u^{(n+1)},v^{(n)})- J(u^{(n+1)},v^{(n+1)}) \geq 0.\notag
\end{align}
Thus, $(J(u^{(n)},v^{(n)}))_n$ is a nonincreasing sequence, 
and since $J \geq 0$ this implies that $(J(u^{(n)},v^{(n)}))_n$ converges. 
Moreover,
$$
J(u^{(0)},v^{(0)}) \geq J(u^{(n)},v^{(n)}) \geq \sum_{\lambda \in \Lambda} \omega_\lambda \|u_\lambda^{(n)}\|_2^2.
$$
Therefore, $(u^{(n)})_n$ is uniformly bounded in
$\ell_{2,\omega^{1/2}}(\Lambda,\mathbb{R}^M)$ and thus,
there exists a subsequence $(u^{(n_k)})_k$ that converges to 
$u^{(\infty)} \in \ell_{2,\omega^{1/2}}(\Lambda,\mathbb{R}^M)$ weakly 
in both $\ell_{2,\omega^{1/2}}(\Lambda,\mathbb{R}^M)$ and $\ell_{2}(\Lambda,\mathbb{R}^M)$, due to our assumption $\omega_\lambda \geq \gamma >0$ for all $\lambda \in \Lambda$. 
For simplicity, let us denote again $u^{(n_k)} = u^{(n)}$. 

First of all, observe that weak convergence implies componentwise convergence, so that $u_\lambda^{(n)} \rightarrow u_\lambda^{(\infty)}$ and $[T^* T u^{(n)}]_\lambda \rightarrow [T^* T u^{(\infty)}]_\lambda$ for all $\lambda \in \Lambda$. By the explicit
formula (\ref{min_v}) for $v_\lambda^{(n)}$ this implies that $v^{(n)}$ 
converges pointwise 
to the limit 
\begin{equation}\label{vlimit}
v^{(\infty)}_\lambda \,:=\, \lim_n v^{(n)}_\lambda
\, = \, \left \{ \begin{array}{ll} \rho_\lambda - \frac{1}{2 \theta_\lambda} \|{u^{(\infty)}}_\lambda\|_q  & \mbox{ if } \|{u^{(\infty)}}_\lambda\|_q  < 2 \theta_\lambda \rho_\lambda,\\
0 & \text{ otherwise}.
\end{array}
\right .
\end{equation}
By definition of $u^{(n)}$ in (\ref{firstalg}) we have $0 \in \partial J_u(u,v^{(n)})$ (where
$\partial J_u(u,v)$ denotes the subdifferential of $J$ considered
as a functional of $u$ only). This means that
$$
0 \in \left[ 2 T^*(T u^{(n)} - g) \right]_\lambda+ 
v_\lambda^{(n-1)} \partial\|\cdot\|_q(u_\lambda^{(n)}) 
+ 2 \omega_\lambda u_\lambda^{(n)} \quad \mbox{ for all } \lambda \in \Lambda,
$$
see also Lemma \ref{lem_subdiffT} and (\ref{subdiff_Phix}), in other words
\begin{equation}
\label{uIncl}
 0 = \left[ 2 T^*(T u^{(n)} - g) \right]_\lambda+ 
v_\lambda^{(n-1)} \zeta_\lambda^{(n)} + 2 \omega_\lambda u_\lambda^{(n)},
\end{equation}
for a suitable $\zeta_\lambda^{(n)} \in \partial\|\cdot\|_q(u_\lambda^{(n)})$.
Now, let 
$(\xi^{(n)}, \eta^{(n)}) \in DJ(u^{(n)},v^{(n)})$.
By definition of $DJ$ and by \eqref{uIncl} we have
\begin{align}
\xi_\lambda^{(n)} \,&\in\, [2T^*(Tu^{(n)}-g)]_\lambda 
+ v_\lambda^{(n)} \partial\|\cdot\|_q(u_\lambda^{(n)}) 
+ 2 \omega_\lambda u_\lambda^{(n)}
\,=\,
v^{(n)}_\lambda \partial \|\cdot\|_q(u_\lambda^{(n)}) - v^{(n-1)}_\lambda \zeta_\lambda^{(n)},\notag
\end{align}
for a suitable $\zeta_\lambda^{(n)} \in \partial\|\cdot\|_q(u_\lambda^{(n)})$.
Since $v^{(n)}_\lambda$ converges 
it is possible to choose the sequence $\xi^{(n)}$ such that
$\lim_{n\to \infty} \xi_\lambda^{(n)} = 0$ for all $\lambda \in \Lambda$.
From (\ref{vlimit}) it is straightforward to check that
\begin{equation}
\label{vIncl}
0 \in \partial s^+(v_\lambda^{(\infty)}) 
\|u_\lambda^{(\infty)}\|_q +2 \theta_\lambda 
(v_\lambda^{(\infty)}-\rho_\lambda) 
\quad \mbox{ for all } \lambda \in \Lambda,
\end{equation}
and similarly 
\[
0 \,\in\, \partial s^+(v_\lambda^{(n)}) \|u_\lambda^{(n)}\|_q
+ 2\theta_\lambda(v_\lambda^{(n)} - \rho_\lambda) 
\]
We can choose $\eta^{(n)}=0$ so that
$\lim_n \eta_\lambda^{(n)} = 0$ for all $\lambda \in \Lambda$. Altogether we conclude that 
$0 \in \lim_n (DJ(u^{(n)},v^{(n)}))_\lambda$ for all $\lambda \in \Lambda$. 
By continuity of $T$ and Lemma \ref{lem_osc} we conclude
\begin{align}
0 &\in \lim_n\left[(2(T^*T(u^{(n)}-g))_\lambda,0) + 
\partial \Phi^{(q)}_\lambda(u_\lambda^{(n)},v_\lambda^{(n)})\right]
\notag\\
& \subset (2 T^*T(u^{(\infty)} -g)_\lambda, 0) + \partial \Phi^{(q)}_\lambda(u_\lambda^{(\infty)},v_\lambda^{(\infty)})
\,=\, DJ(u^{(\infty)},v^{(\infty)})_\lambda\notag
\end{align}
for all $\lambda \in \Lambda$.
It follows that $0 \in DJ(u^{(\infty)},v^{(\infty)}) \subset
\partial J(u^{(\infty)},v^{(\infty)})$, the latter inclusion by Proposition (\ref{subdiff}). 
Hence, by strict convexity $(u^*,v^*) = (u^{(\infty)},v^{(\infty)})$.
With this we have shown the weak convergence of the sequence $u^{(n)}$ to $u^*$. 

To establish the strong convergence we need to develop a more detailed analysis of the minimization of $J$ with respect to $u$. Next section is devoted to this end, and it will allow us to use some further tools for the full proof of Theorem \ref{thm_conv} in Subsection 5.3.

\section{An Iterative Thresholding Algorithm for the Minimization with Respect to $u$}
\label{sec_iter}

One step of the minimization algorithm in the previous section consists 
in minimizing $J(u,v)=J^{(q)}_{\theta,\rho,\omega}(u,v)$ for some fixed $v$. Moreover, keeping 
$v$ fixed is also
interesting for its own -- in particular, if one is interested 
in minimizing the functional $K = K^{(q)}_v$ defined in (\ref{def_Psi}).
Indeed, for $\omega=0$ and $\rho = v$ 
we have $J^{(q)}_{\theta,v,0}(u,v) = K^{(q)}_v(u)$. 
As we will describe in the following this minimization task
can be performed by a thresholded Landweber algorithm similar
to the one analyzed by Daubechies et al. in \cite{dadede04}.

With $v$ fixed our task is equivalent to minimizing
\begin{equation}\label{def_K}
K(u) \,=\, K^{(q)}_{v,\omega}\,:=\, \|T u - g\|_\cH^2 + \Psi(u)
\end{equation}
with respect to $u \in \ell_2(\Lambda,\RR^M)$ where
\begin{equation}\label{def_Psi}
\Psi(u) \,:=\, \Psi^{(q)}_{v,\omega}(u) \,:=\, 
\sum_{\lambda \in \Lambda} v_\lambda \|u_\lambda\|_q + 
\sum_{\lambda \in \Lambda} \omega_\lambda \|u_\lambda\|_2^2.
\end{equation}
We assume that $T$ is non-expansive, i.e., $\| T \| < 1$, which can always 
be achieved by rescaling. Also we suppose that $K$ is strictly convex.
This is ensured if e.g. the kernel of $T$ is trivial or $\omega_\lambda > 0$
for all $\lambda \geq 0$. 

We define a surrogate functional by
\[
K^{s}(u,a) \,:=\, K(u) - \|Tu - Ta\|^2_\cH + \|u-a\|^2_2
\,=\, \|Tu-g\|^2_\cH + \Psi(u) - \|Tu - Ta\|^2_\cH + \|u-a\|^2_2.
\]
Since $\|T\| < 1$ also $K^s$ is convex, see \cite{dadede04} for a rigorous
argument. Now starting with some $u^{(0)} \in \ell_2(\Lambda,\RR^M)$
we define a sequence $u^{(m)}$ by
\[
u^{(m+1)} = \argmin_{u \in \ell^2(\Lambda,\RR^M)} K^s(u,u^{(m)})
\]
The minimizer of $K^s(u,a)$ (for fixed $a$) can be determined
explicitly as follows. First, we claim that
\[
\argmin_u K^s(u,a) \,=\, U_\Psi(a + T^*(g- T a)),
\]
where the ``thresholding'' operator $U_\Psi$ is defined as
\begin{equation}\label{def_UPsi}
U_\Psi(u) \,:=\, \argmin_{z \in \ell_2(\Lambda,\RR^M)} \|u-z\|^2_2 + \Psi(z).
\end{equation}  
Indeed, a direct calculation shows that
\begin{align}
K^s(u,a) \,&=\, \|Tu - g\|^2_\cH - \|Tu-Ta\|^2_\cH + \|u-a\|^2_2 + \Psi(u)\notag\\
&=\, \|(a+T^*(g-Ta))-u\|^2_2 + \Psi(u) - \|a+T^*(g-Ta)\|^2_2+\|g\|^2_\cH-\|Ta\|^2_\cH + \|a\|^2_2.\notag
\end{align}
Since the last terms (after $\Psi(u)$) do not depend on $u$ they can be discarded when
minimizing with respect to $u$, and the above claim follows. (The same argument works
also for general 'sparseness measures' $\Psi$).
Thus, the iterative algorithm reads
\begin{equation}\label{def_un}
u^{(m+1)} \,=\, U_\Psi(u^{(m)} + T^*(g-T u^{(m)})).
\end{equation}
In the following we give more details about $U_\Psi$ and analyze the
convergence of this algorithm.

\subsection{The thresholding operator}

Let us derive more information about $U_\Psi$ for our specific $\Psi=\Psi_{v,\omega}$
in (\ref{def_Psi}). We have the following lemma.
\begin{lemma}\label{lem_thresh1} Let $1\leq q \leq \infty$. It holds
\[
(U_{v,\omega}^{(q)}(u))_\lambda \,:=\, (U_{\Psi^{(q)}_{v,\omega}} (u))_\lambda 
\,=\, (1+\omega_\lambda)^{-1} S^{(q)}_{v_\lambda}(u_\lambda),
\] 
where
\begin{equation}\label{def_Svq}
S^{(q)}_{v}(x) \,=\, \argmin_{z\in \RR^M} \|z-x\|_2^2 + v\|z\|_q,\quad x\in \RR^M.
\end{equation}
Furthermore, $S^{(q)}_v$ is given by
\begin{equation}\label{Svq_explicit}
S^{(q)}_v(x) \,=\, x - P^{q'}_{v/2}(x),
\end{equation}
where $P^{q'}_{v/2}$ denotes the orthogonal projection onto the norm ball of radius $v/2$
with respect to the dual norm of $\|\cdot\|_{q}$, i.e., the $\|\cdot\|_{q'}$-norm 
with $q'$ denoting the dual index, $1/q+1/q' = 1$.  
(The analogous result holds also if the norm $\|\cdot\|_q$ is replaced by an arbitrary norm on $\RR^M$).
\end{lemma}
\begin{proof} For $\Psi_{v,\omega}$ the minimizing problem defining $U_\Psi=U^{(q)}_{v,\omega}$ 
decouples with respect to $\lambda \in \Lambda$. Thus, we have
\[
(U^{(q)}_{v,\omega}(x))_\lambda = \argmin_{z \in \RR^M} \|x_\lambda-z\|_2^2 + \omega_\lambda \|z\|_2^2+ v_\lambda\|z\|_q.
\]
If $z$ minimizes the latter term then necessarily
$
0 \in 2(1+\omega_\lambda)z- 2x + v_\lambda \partial \|\cdot\|_q(z) 
$
where $\partial\|\cdot\|_q$ denotes the subdifferential of the $q$-norm. In other words, 
\[
(1+\omega_\lambda)z - x \in - \frac{v_\lambda}{2} \partial \|\cdot\|_q(z).
\]
Since $\|\cdot\|_q$ is $1$-homogeneous we have $\partial \|\cdot\|_q(z) = \partial \|\cdot\|_q((1+\omega_\lambda)z)$.
Setting $y=(1+\omega_\lambda)z$ gives $y-x \in -\frac{v_\lambda}{2} \partial \|\cdot\|_q(y)$, which is the above
relation for $\omega_\lambda = 0$. From this we deduce the first claim.

Let us show the second claim, i.e., the explicit form of the operator $S_v^{(q)}$. We already know that if $z$
minimizes the left hand side of (\ref{def_Svq}) then $x-z \in \partial \frac{v}{2}\|z\|_q$. Let
$\psi(z) = \frac{v}{2}\|z\|_q$ and $\psi^*$ be its Fenchel conjugate function defined by
$\psi^*(y) \,=\, \sup_{x} (\langle x,y \rangle - f(x))$. It is well-known \cite[p.~93]{bova04} that
\[
\psi^*(y) \,=\, \chi_{B^{q'}(v/2)}(y) \,:=\, 
\left\{ \begin{array}{ll} 0 & \mbox{if } \|y\|_{q'} \leq v/2 \\
\infty & \mbox{otherwise} \end{array} \right.
\]
Here $B^{q'}(v/2)$ denotes the norm ball of radius $v/2$ with respect to 
the dual norm of $\|\cdot\|_q$. 
It is a standard result, see e.g. \cite[Proposition 11.3]{rowe98}, \cite[Corollary 5.2]{ekte99}, that
$w \in \partial \psi(y)$ if and only if $y \in \partial \psi^*(w)$ yielding
$z \in \partial \psi^*(x-z)$ in our case, and hence,
\[
x \in x-z + \partial \psi^*(x-z) \,=\, x-z +\partial \chi_{B^{q'}(v/2)}(x-z).
\]
Now if $y \in w + \partial \chi_{B^{q'}(v/2)}(w)$ then it is straightforward to see
that $w$ must be the orthogonal projection of $y$ onto $B^{q'}(v/2)$, i.e.,
$w = \argmin_{w' \in B^{q'}(v/2)} \|w'-y\|_2$, see also \cite[Example 10.2 and p.~20]{rowe98}. 
For our situation this means that
$x-z = P^{q'}_{v/2}(x)$, i.e., $z = x - P^{q'}_{v/2}(x)$. This shows the second claim.

Clearly, all arguments work also for a general norm rather than the $q$-norm.
\end{proof}
Let us give $S^{(q)}_v$ explicitly for $q=1,2,\infty$.

\begin{lemma}\label{explS} Let $x \in \RR^M$ and $v \geq 0$.
\begin{itemize}
\item[(a)] For $q=1$ we have $S^{(1)}_v(x) = (s^{(1)}_v(x_\ell))_{\ell = 1}^M$ where for $y \in \RR$
\[
s^{(1)}_{v}(y) \,=\,  \left\{\begin{array}{ll} 0 & \mbox{ if } |y| \leq \frac{v}{2},\\
\sign(y)(|y|-\frac{v}{2}) & \mbox{ otherwise.}\end{array} \right.
\]
\item[(b)] For $q=2$ it holds
\[
 S^{(2)}_{v}(x) :=\left\{\begin{array}{ll} 0 & \mbox{ if }  \|x\|_2 \leq \frac{v}{2},\\
\frac{(\|x\|_2 - v/2)}{\|x\|_2} x &  \mbox{ otherwise.} \end{array}\right.
\]
\item[(c)] Let $q=\infty$.
Order the entries of $x$ by magnitude
such that $|x_{i_1}| \geq |x_{i_2}| \geq \hdots \geq |x_{i_M}|$. 
\begin{enumerate}
\item If $\|x\|_1 < v/2$ then $S_v^{(\infty)}(x) = 0$.
\item If $\|x\|_1 > v/2$, 
let $n \in \{1,\hdots,M\}$ be the largest index satisfying
\begin{equation}\label{cond_infty}
|x_{i_n}| \,\geq\, \frac{1}{n-1}\left(\sum_{k=1}^{n-1} |x_{i_k}| - \frac{v}{2}\right).
\end{equation}
Then
\begin{align}
(S_v^{(\infty)}(x))_{i_j} \,&=\, 
\frac{\sign(x_{i_j})}{n}\left(\sum_{k=1}^n |x_{i_k}| - \frac{v}{2}\right), \quad j=1,\hdots,n, \notag\\
(S_v^{(\infty)}(x))_{i_j} \,&=\, x_{i_j},\quad j=n+1,\hdots,M.\notag
\end{align}
\end{enumerate}
\end{itemize}
\end{lemma}
\begin{proof} (b) The projection $P^2_{v/2}(x)$ of $x$ onto an $\ell_2$ ball of radius $v/2$ is
clearly given by
\[
P^2_{v/2}(x) \,=\, \left\{ \begin{array}{ll} x & \mbox{ if } \|x\|_2 \leq v/2,\\
\frac{v/2}{\|x\|_2}x & \mbox{ otherwise }.
\end{array}\right.
\] 
Since by the previous lemma $S_v^{(2)}(x) = x - P^2_{v/2}(x)$ this gives the assertion.

(a) Although this is well-known we give a simple argument. For $q=1$ the functional in (\ref{def_Svq})
defining $S^{(1)}$ decouples, i.e.,
\[
S_v^{(1)}(x) \,=\, \argmin_{z \in \RR^M} \sum_{\ell=1}^M \left(|z_\ell - x_\ell|^2 + v|z_\ell|\right).
\]
Thus, $S_v^{(1)}(x)_\ell = \argmin_{z_\ell \in \RR} |z_\ell - x_\ell|^2 + v|x_\ell|$ for all $\ell=1,\hdots,M$. 
The latter can be interpreted as the problem for $q=2$ on $\RR^1$ and hence, the assertion follows from (b).

(c) 
If $\|x\|_1 \leq v/2$ then $P^1_{v/2}(x) = x$ and by the previous lemma 
$S^{(\infty)}_v(x) = x - P^1_{v/2}(x) = 0$. 
Now assume $\|x\|_1 > v/2$. Let $z = S^{(\infty)}_{v}(x)$. This is equivalent to
$0$ being contained in the subdifferential of the functional in (\ref{def_Svq}) defining
$S^{(\infty)}_v$. This means
\begin{equation}\label{x_rel}
2(z-x) \in - v \partial \|\cdot\|_\infty(z).
\end{equation} 
We recall that the subdifferential of the maximum norm is given by (\ref{max_subgrad}).

Now assume for the moment
that the maximum norm of $z$ is attained in $z_{i_1},\hdots,z_{i_n}$.
We will later check whether this was really the case. 
Further, we assume for simplicity
that all the entries $x_{i_1},\hdots,x_{i_n}$ are positive. (The other cases
can be carried through in the same way).
Then certainly also the numbers $z_{i_1},\hdots,z_{i_n}$ are positive because
choosing them with the opposite sign would certainly increase the functional defining $S^{(\infty)}_v$. 
Then by (\ref{max_subgrad}) we obtain
$2(z_{i_j} - x_{i_j}) = 0$
for the entries $z_{i_j}$ not giving the maximum, i.e.,
$
z_{i_j} \,=\, x_{i_j}, j=n+1,\hdots,M.
$

Moreover, if $n=1$ (i.e., the maximum norm
of $z$ is attained at only one entry) then
$2(z_{i_1} - x_{i_1}) = -v$, in other words,
$
z_{i_1} \,=\, x_{i_1} - v/2.
$
Thus, the initial hypothesis
that the maximum norm of $z$ is attained only at $z_{i_1}$ is true if and only if 
the second largest entry
$x_{i_2}$ satisfies $|x_{i_2}| < |z_{i_1}| - v/2$.

So if the latter inequality is not satisfied then the maximum norm of $z$ is at least attained at
two entries, i.e., $n\geq 2$. In this case by (\ref{max_subgrad}) 
the entries $z_{i_1}=z_{i_2} = \cdots = z_{i_n} = t$ 
satisfy
\begin{align}
2 t - 2 x_{i_j} \,& =\, - v a_j,\quad j=1,\hdots,n-1,\notag\\
2 t - 2 x_{i_n} \,& =\, - v \left(1-\sum_{k=1}^{n-1} a_k\right) \notag
\end{align}
for some numbers $a_1,\hdots,a_{n-1} \in [0,1]$ satisfying $\sum_j a_j \leq 1$. 
This is a system of $n$ linear equations
in $t$ and $a_1,\hdots,a_{n-1}$. Writing it in matrix form we get
\[
\left(\begin{matrix}
1 & v/2 & 0 & 0 &\cdots & 0\\
1 & 0  & v/2 & 0 & \cdots & 0\\
\vdots & \vdots &\vdots & \vdots & \vdots & \vdots\\
1 & -v/2 & -v/2 & - v/2 & \cdots & -v/2
\end{matrix} \right)
\left(\begin{matrix} t\\ a_1 \\  \vdots \\ a_{n-1}\end{matrix}\right)
\,=\, \left(\begin{matrix}
x_{i_1}\\ \vdots \\ x_{i_{n-1}} \\ x_{i_n} - v/2 \end{matrix} \right).
\]
Denoting the matrix on the left hand side by $B$, a simple computation verifies that
\[
B^{-1} \,=\, \frac{1}{n} \left(\begin{matrix}
1 & 1 & 1 & \cdots & 1\\
\frac{2(n-1)}{v} & -\frac{2}{v} & - \frac{2}{v} & \cdots & - \frac{2}{v} \\
-\frac{2}{v} & \frac{2(n-1)}{v} & - \frac{2}{v} & \cdots & - \frac{2}{v} \\
\vdots & \vdots & \ddots & \vdots & \vdots \\
-\frac{2}{v} & \cdots & -\frac{2}{v} & \frac{2(n-1)}{v} & - \frac{2}{v}
\end{matrix}\right).
\]
This gives
\[
z_{i_1} \,=\, \hdots \,=\, z_{i_n} \,=\, t \,=\,
\frac{1}{n}\left(\sum_{j=1}^n x_{i_j} - v/2\right)
\] 
and
$
a_j = \frac{2}{nv}\left(v/2 + (n-1) x_{i_j} - 
\sum_{k \in \{1,\hdots,n\} \setminus \{j\}} x_{i_k}\right).
$
Thus, all $a_j$ are non-negative if for all $j \in \{1,\hdots,n-1\}$
\[
x_{i_j} \geq \frac{1}{n-1}
\left( \sum_{k \in \{1,\hdots,n\} \setminus \{j\}} x_{i_k} -v/2\right).
\]
Moreover, a simple calculation gives
$
\sum_{j=1}^{n-1} a_j = \frac{n-1}{n} + \frac{2}{nv}\left(\sum_{j=1}^{n-1} x_{i_j}
-(n-1) x_{i_n}\right).
$
Thus, it holds $1-\sum_{j=1}^{n-1} a_j \geq 0$ if and only if
\[
x_{i_n} \,\geq\, \frac{1}{n-1}\left(\sum_{j=1}^{n-1} x_{i_j} - v/2\right).
\]
Therefore, the initial assumption that the maximum norm of $u$ is attained precisely at 
$z_{i_1},\hdots,z_{i_n}$ can only be true if $x_{i_1},\hdots,x_{i_n}$
are the largest entries of the vector $x$ and 
\[
z_{i_{n+1}} \,=\, x_{i_{n+1}} \,<\, t \,=\, 
\frac{1}{n} \left(\sum_{j=1}^n x_{i_j} - v/2\right),
\]
i.e., $|x_{i_{n+1}}| < n^{-1}(\sum_{j=1}^n |x_{i_j}| - v/2)$.
Pasting all the pieces together shows the assertion of the lemma.
\end{proof}

\subsection{Weak convergence}

In the following we will prove that $u^{(m)}$ converges weakly and 
strongly to the unique minimizer of $K$.
We first establish the weak convergence. Following the proof of Proposition 3.11 in \cite{dadede04}
one may extract essentially three conditions on a general sparsity measure $\Psi$ such that
weak convergence is ensured. Let us collect them in the following Proposition.

\begin{proposition}\label{prop_weak} Assume $K$ is given by (\ref{def_K}) with a general sparsity measure $\Psi$ and suppose $K$ is strictly convex.
Let $U_\Psi$ be the associated 'thresholding operator' given by (\ref{def_UPsi}). Assume that the following
conditions hold
\begin{itemize}
\item[(1)] $U_\Psi$ is non-expansive, i.e. $\|U_\psi(x) - U_\Psi(y)\|_2 \leq \|x-y\|_2$ for all 
$x,y \in \ell_2(\Lambda,\RR^M)$.
\item[(2)] It holds $\|f\|_2 \leq H(\Psi(f))$ for all $f \in \ell_2(\Lambda,\RR^M)$ and some 
monotonically increasing function $H$ on $\RR_+$.
(This ensures that a sequence $f_n$ satisfying $\Psi(f_n) \leq  C$ is bounded in $\ell_2(\Lambda,\RR^M)$).
\item[(3)] For all $x,h \in \ell_2(\Lambda,\RR^M)$ it holds
\[
\Psi(U_\Psi(x)+h) - \Psi(U_\Psi(x)) + 2 \langle h, U_\Psi(x) -x \rangle \geq 0.
\] 
\end{itemize}
Then the sequence $u^{(m)}$ defined by (\ref{def_un}) converges
weakly to the minimizer of $K$ independently of the choice of $u^{(0)}$.
\end{proposition}
\begin{proof} First we claim that the condition in (1) implies that the surrogate
functional $K^s$ satisfies
\begin{equation}\label{rel_surrogate}
K^s(u+h,a) - K^s(u,a) \,\geq\, \|h\|_2^2
\end{equation}
for $u = \argmin_{u'} K^s(u',a)\,=\, U_\Psi(a-T^*(g-Ta))$. Indeed, set $x := a-T^*(g-Ta)$, i.e.,
$u = U_\Psi(x)$. Then an elementary calculation yields
\begin{align}
K^s(u+h,a) - K^s(u)
=\, & \|T(u+h) - g\|_2^2 + \Psi(u+h) - \|T(u+h)-Ta\|_2^2 + \|u+h-a\|_2^2\notag\\
& -  \|Tu-g\|_2^2 - \Psi(u) + \|Tu-Ta\|_2^2 - \|u-a\|_2^2\notag\\
=\, & 2\langle h, u-a-T^*(g-Ta)\rangle + \Psi(u+h) - \Psi(u) + \|h\|^2_2\notag\\
=\, & 2\langle h, U_\Psi(x) -x \rangle + \Psi(U_\Psi(x) + h) - \Psi(U_\Psi(x)) + \|h\|_2^2 \,\geq\, \|h\|_2^2.\notag
\end{align}
The relation in (1) was used in the last inequality. 

Now with (\ref{rel_surrogate}) and properties (2) and (3) one can easily justify
that the proofs of the analogues of Theorem 3.2 until Proposition 3.11 in 
\cite{dadede04} go through completely in the same way, which finally leads 
to the statement of this proposition.
\end{proof}

Let us now show that for our specific choice of $\Psi = \Psi^{(q)}_{v,\omega}$ properties
(1) - (3) in the previous Proposition hold, and thus, $u^{(n)}$ converges weakly to
a minimizer of $K$.

\begin{lemma}\label{nonexp} $U_{v,\omega}^{(q)}=U_{\Psi^{(q)}_{v,\omega}}$ is non-expansive.
\end{lemma}
\begin{proof} Clearly, the map $x \mapsto (1+\omega_\lambda)^{-1} x$ is non-expansive. 
By (\ref{Svq_explicit}) we have
$S^{(q)}_v = I - P_{v/2}^{q'}$. Since $P_{v/2}^{q'}$ is an orthogonal projection onto a convex
set also $S^{(q)}_v$ is non-expansive, see e.g. \cite{te05}. Hence, $U^{(q)}_{v,\omega}$ 
is non-expansive
since on each component $x_\lambda$, $\lambda \in \Lambda$, it is a 
composition of non-expansive operators.
\end{proof}

\begin{lemma} If $(v_\lambda)$ or $(\omega_\lambda)$ are bounded away from $0$ then 
condition (2) in Proposition \ref{prop_weak} holds.
\end{lemma}
\begin{proof} This follows by a standard argument.
\end{proof}

If we consider the problem of minimizing $J(u,v)$ jointly over $u$ and $v$ then we certainly
cannot assume that $v$ is bounded away from $0$, but in this case we require
that $\omega_\lambda$ is bounded away from $0$. (By Proposition
\ref{q12oo} this is needed anyway to ensure that $J(u,v)$ is jointly convex in
$u$ and $v$).
In the case where we only
minimize $J(u,v)$ with respect to $u$ (i.e., when minimizing $\Psi(u)$ defined in (\ref{def_Psi})) 
we may take $\omega_\lambda$ arbitrary
(and even $\omega_\lambda = 0$) but then we have to require a lower
bound on $v_\lambda$.

Now consider the third condition in the Proposition. The next lemma shows
that it suffices to prove it for $S_v^{(q)}$, i.e., for $\omega_\lambda =0$.

\begin{lemma} Assume that for all $x,h \in \RR^M$ it holds
\[
v(\|S_v^{(q)}(x) + h\|_q - \|S_v^{(q)}(x)\|_q) +2 \langle h, S_v^{(q)}(x) - x \rangle \geq 0.
\] 
Then condition (3) in Proposition \ref{prop_weak} is satisfied.
\end{lemma}
\begin{proof} By definition of $\Psi_{v,\omega}^{(q)}$ we need to show that for $\omega, v \geq 0$
and all $x,h \in \RR^M$
\begin{align}
&v (\|(1+\omega)^{-1} S_{v}^{(q)}(x) +h \|_q - \|(1+\omega)^{-1} S_v^{(q)}(x) \|_q)\notag\\
+ & \omega(\|(1+\omega)^{-1}S_v^{(q)}(x)+h\|_2^2 - \|(1+\omega)^{-1} S_v^{(q)}(x)\|_2^2)
+ 2\langle h, (1+\omega)^{-1} S_v^{(q)}(x) -x\rangle \geq 0.   \notag
\end{align}
Setting $h' = (1+\omega)h$ we obtain for the left hand side of this inequality
\begin{align}
&(1+\omega)^{-1}v\left(\|S_v^{(q)}(x) + h'\|_q - \|S_v^{(q)}(x)\|_q\right)\notag\\
+ & (1+\omega)^{-2} \omega\left(\|S_v^{(q)}(x) + h'\|_2^2 - \|S_v^{(q)}(x)\|_2^2\right)
+ 2(1+\omega)^{-2} \langle h', S_v^{(q)}(x)-x\rangle\notag\\
= & (1+\omega)^{-1} \left[v(\|S_v^{(q)}(x) + h'\|_q - \|S_v^{(q)}(x)\|_q) +2\langle h', S_v^{(q)}(x)-x\rangle\right]\notag\\
+ & (1+\omega)^{-2} \omega\left[\|S_v^{(q)}(x) + h'\|_2^2 - \|S_v^{(q)}(x)\|_2^2 
- 2\langle h', S_v^{(q)}(x) \rangle\right] \geq (1+\omega)^{-2} \omega \|h'\|_2^2 \geq 0.   \notag
\end{align}
This completes the proof.
\end{proof} 


\begin{lemma} The condition in the previous lemma holds for $S_v^{(q)}$, $1\leq q \leq \infty$
(and even if the $\ell_q$ norm is replaced by a general norm on $\RR^M$).
\end{lemma}
\begin{proof} First note that by definition (\ref{def_Svq}) and duality we have
\[
\|S_v^{(q)}(x)\|_q \,=\, (v/2)^{-1} \sup_{k \in B^{q'}(v/2)} \langle k, x-P_{v/2}^{q'} x\rangle
\]
A characterization of the orthogonal projection 
tells us that
$\langle k-P_{v/2}^{q'} x,x-P_{v/2}^{q'} x\rangle \leq 0$ for all $k \in B^{q'}(v/2)$, see e.g.~\cite[Lemma 8]{te05}.
This gives
\begin{align}
 \|S_v^{(q)}(x)\|_q \,&=\, (v/2)^{-1} \sup_{k \in B^{q'}(v/2)} \left(\langle k-P_{v/2}^{q'} x, x-P_{v/2}^{q'} (x)\rangle 
+ \langle P_{v/2}^{q'}(x), S_v^{(q)}(x)\rangle\right)\notag\\
 &\leq\,
(v/2)^{-1} \langle P_{v/2}^{q'}(x), S_{v}^{(q)}(x)\rangle. \notag
\end{align} 
Using once more that $S_v^{(q)}(x) - x = - P_{v/2}^q(x)$ we further obtain
\begin{align}
&\,v(\|S_v^{(q)}(x)+h\|_q - \|S_v{(q)}(x)\|_q) + 2\langle h, S_v^{(q)}(x)-x \rangle \notag\\
\geq & \, v\|S_v^{(q)}(x)+h\|_q - 2\langle P_{v/2}^{q'}(x), S_v^{(q)}(x) \rangle - 2\langle P_{v/2}^{q'}(x),h\rangle 
= v \|S_v^{(q)}(x)+h\|_q - 2\langle P_{v/2}^{q'}(x), S_v^{(q)}(x) + h \rangle\notag\\
\geq &\, v \|S_v^{(q)}(x)+h\|_q - 2\|P^{q'}_{v/2}x\|_{q'} \|S_v^{(q)}(x) + h\|_q
\,\geq\, 0.  \notag
\end{align}
Hereby, we used that $P^{q'}_{v/2}$ is a projection onto $B^{q'}(v/2)$, so $\|P^{q'}_{v/2}x\|_{q'} \leq v/2$.
This finishes the proof.
\end{proof}

To summarize we have the following result about weak convergence.

\begin{corollary} Let $1\leq q \leq \infty$ and assume that $(v_\lambda)$ or $(\omega_\lambda)$ is bounded
from below. Then the sequence $u^{(m)}$ defined in (\ref{def_un}) converges weakly to a minimizer
of $K$, where $\Psi= \Psi_{v,\omega}^{(q)}$ is the sparsity measure defined in (\ref{def_Psi}). 
(The $q$-norm in (\ref{def_Psi}) can be replaced by any other norm on $\RR^M$).
\end{corollary}

\subsection{Strong convergence}

The next result establishes the strong convergence.

\begin{proposition} Let $1\leq q \leq \infty$ and assume that $(v_\lambda)$ or $(\omega_\lambda)$ are bounded
away from $0$. In case $(v_\lambda)$ is not bounded away from $0$ assume further that 
there is a constant $c>0$
such that $v_\lambda < c$ for only finitely many $\lambda$.
Then $u^{(m)}$ converges strongly to a minimizer of $K$.
\end{proposition}
\begin{proof} The analogues of Lemmas 3.15 and 3.17 in \cite{dadede04}
are proven in completely the same way. It remains to justify the analogue
of \cite[Lemma 3.18]{dadede04}: If for some $a \in \ell_2(\Lambda,\RR^M)$ and some sequence
$(h^{(m)}) \subset \ell_{2}(\Lambda,\RR^M)$ converging weakly to $0$ it holds 
$\lim_{m\to \infty}
\|U^{(q)}_{v,\omega}(a + h^{(m)}) - U^{(q)}_{v,\omega}(a) - h^{(m)}\|_2 = 0$
then $\|h^{(m)}\|_2 \to 0$ for $m\to \infty$. To this end we mainly follow
the argument in \cite{dadede04}.

Let $c$ be the constant such that $v_\lambda < c$ for 
$\lambda \in \Lambda_{00}$ for $\Lambda_{00}$ finite. Then let
$\Lambda_{01}$ be a finite set such that 
$\sum_{\lambda \in \Lambda \setminus \Lambda_{01}} 
\|a_\lambda\|_{q'} \leq \sigma$ for some $\sigma < c/2$.
(Such a set $\Lambda_{01}$ exists since $\|\cdot\|_{q'}$ and $\|\cdot\|_2$ are equivalent
norms on $\RR^M$ and by assumption $a \in \ell_2(\Lambda,\RR^M)$). Since
$\Lambda_0 = \Lambda_{00} \cup \Lambda_{01}$ is also finite,
we have $\sum_{\lambda \in \Lambda_0} \|h^{(m)}_\lambda\|_2^2 \to 0$
for $m\to \infty$ by the weak convergence of $h^{(m)}$ to $0$. Thus,
we are left with proving that $\sum_{\lambda \in \Lambda \setminus \Lambda_0} \|h^{(m)}_\lambda\|_2^2 \to 0$
for $m \to \infty$.

For each $m$ we split $\Lambda_1 := \Lambda \setminus \Lambda_0$ into the subsets
$\Lambda_{1,m} := \{\lambda \in \Lambda_1: \|h^{(m)}_\lambda + a_\lambda\|_{q'} < v_\lambda/2\}$
and $\widetilde{\Lambda}_{1,m} = \Lambda_1 \setminus \Lambda_{1,m}$. If $\lambda \in \Lambda_1$
then $U_{v,\omega}^{(q)}(a + h^{(m)})_\lambda = U^{(q)}_{v,\omega}(a)_\lambda = 0$
since $\|a_\lambda + h^{(m)}_\lambda\|_{q'}, \|a_\lambda\|_{q'} \leq v_\lambda/2$.
Thus, $\|h^{(m)}_\lambda - U_{v,\omega}^{(q)}(a + h^{(m)})_\lambda + U^{(q)}_{v,\omega}(a)_\lambda\|_2^2
= \|h^{(m)}_\lambda\|_2^2$ and by assumption,
\[
\sum_{\lambda \in \Lambda_1} \|h^{(m)}_\lambda\|_2^2 
\leq \|h^{(m)} - U^{(q)}_{v,\omega}(a+h^{(m)}) + U^{(q)}_{v,\omega}(a)\|_2^2 \to 0
\quad \mbox{as }m \to \infty.
\]
Now let $\lambda \in \widetilde{\Lambda}_{1,m}$. We first consider the case 
that $\omega_\lambda = 0$, i.e., $U_{v,\omega}{q}(x)_\lambda = S_v^{(q)}(x)_\lambda$. 
Since $\|a_\lambda\|_{q'} \leq \sigma < v_\lambda/2$ we have $S_v^{(q)}(a) = 0$, and thus,
\begin{align}
&\| h^{(m)}_\lambda - S_{v_\lambda}^{(q)}(h^{(m)}_\lambda+a_\lambda)\|_{q'} \,=\, 
\|h^{(m)}_\lambda - (h^{(m)}_\lambda+a_\lambda) + P_{v_\lambda/2}^{q'}(h^{(m)}_\lambda + a_\lambda)\|_{q'}\notag\\
=\,& \|P_{v_\lambda/2}^{q'}(h^{(m)}_\lambda + a_\lambda) - a_\lambda\|_{q'}
\,\geq\, \|P_{v_\lambda/2}^{q'}(h^{(m)}_\lambda+a_\lambda)\|_{q'} - \|a_\lambda\|_{q'}
\,\geq\, v_\lambda/2 - \sigma \geq c/2 - \sigma.\notag
\end{align}
Hereby, we used that $\|P_{v_\lambda/2}^{q'}(h^{(m)}_\lambda+a_\lambda)\|_{q'} = v_\lambda/2$ (because 
$\|h^{(m)}_\lambda+a_\lambda\|_{q'} \geq v_\lambda/2$). Since every norm on a finite-dimensional
space is equivalent there is a constant $C$ such that
\[
\| h^{(m)}_\lambda - S_{v_\lambda}^{(q)}(h^{(m)}_\lambda+a_\lambda)+S_{v_\lambda}^{(q)}(a_\lambda)\|_2^2 \geq C^2(c/2-\sigma)^2 > 0.
\]
However, since by assumption
$\sum_{\lambda \in \Lambda} \| h^{(m)}_\lambda - S_{v_\lambda}^{(q)}(h^{(m)}_\lambda+a_\lambda)+S_{v_\lambda}^{(q)}(a)\|_2^2 \to 0$
as $m \to \infty$ there must exist an $m_0$ such that $\widetilde{\Lambda}_{1,m}$ 
is empty for all $m \geq m_0$.

In the case that $\omega_\lambda$ does not vanish we have
\begin{align}
\| h^{(m)}_\lambda - U^{(q)}_{v,\omega}(h^{(m)} + a)_\lambda +
U^{(q)}_{v,\omega}(a)_\lambda\|_2
\,&=\, \|h^{(m)}_\lambda - 
(1+\omega_\lambda)^{-1}S_{v_\lambda}^{(q)}(h^{(m)}_\lambda + a_\lambda)\|_2\notag\\
&=\, (1+\omega_\lambda)^{-1}\|(1+\omega_\lambda)h^{(m)}_\lambda - S_{v_\lambda}^{(q)}
(h^{(m)}_\lambda + a_\lambda)\|_2
\end{align}
We claim that
\begin{equation}\label{claim0}
\|(1+\omega_\lambda)h^{(m)}_\lambda - S_{v_\lambda}^{(q)}
(h^{(m)}_\lambda + a_\lambda)\|_2 \,\geq\, \|h^{(m)}_\lambda - S_{v_\lambda}^{(q)}(h^{(m)}_\lambda+a_\lambda)\|_2
\end{equation}
so that we can apply the argument for $\omega_\lambda = 0$ 
to conclude that $\widetilde{\Lambda}_{1,m}$
is empty for $m$ sufficiently large.
Let us omit for the moment all indexes $\lambda$ and $m$ 
for the sake of simpler notation. We have
\begin{equation}\label{claim1}
\|(1+\omega)h - S_v^{(q)}(h+a)\|_2^2 - \|h - S_v^{(q)}(h+a)\|_2^2
\,=\, 2\omega\langle h, h-S_v^{(q)}(h+a)\rangle + \omega^2\|h\|_2^2
\end{equation}
and furthermore,
\begin{align}
\langle h, h-S_v^{(q)}(h+a)\rangle \,&=\, \langle h, P_{v/2}^{q'}(h+a)-a\rangle\notag\\
&=\, -\langle h + a - P_{v/2}^{q'}(h+a),a-P_{v/2}^{q'}(h+a)\rangle
+ \|a-P_{v/2}^{q'}(h+a)\|_2^2
\geq 0.
\end{align}
Hereby, we used that $a \in B^{q'}(\sigma) \subset B^{q'}(v/2)$
and the fact that $\langle k - P^{q'}_{v/2}(x),x-P_{v/2}^{q'}(x)\rangle\leq 0$ 
for all
$k \in B^{q'}(v/2)$ and $x \in \RR^M$. Thus, the term in (\ref{claim1})
is non-negative and therefore our claim (\ref{claim0}) holds.
\end{proof}

Let us shortly comment on the condition that if $v_\lambda$ is not bounded from below
there is at least some $c > 0$ such that $v_\lambda > c$ except for
a finite set of indexes $\lambda$. This condition
is mainly relevant when considering also a minimization over $(v_\lambda)$.
Then the term $\sum \theta_\lambda (\rho_\lambda - v_\lambda)^2$ in the functional
$J(u,v)$ ensures that the sequence $(\rho_\lambda - v_\lambda)$ is contained in
$\ell_{2,\theta^{1/2}}$. If $\theta_\lambda$ and $\rho_\lambda$ are bounded from below this implies
that $v_\lambda$ can be less than $1/2 \min_\lambda \rho_\lambda$, say, only for finitely
many $\lambda$.

\section{Numerical Implementation and Error Analysis}
\label{sec_image}

The scope of this section is twofold: We want to formulate an implementable version 
of the double-minimization algorithm and show its strong convergence. To this end we 
develop an error analysis.

\subsection{Numerical implementation}

Let us compose the two iterative algorithms described 
in \eqref{firstalg} and \eqref{def_un}, respectively, into a unique scheme.
 
\begin{alg}
  \label{alg1}
  \noindent {\em
    {\bf JOINTSPARSE}\\ 
    \begin{tabular}{ll}
    Input:    & Data vector $(g_j)_{j=1}^N$, initial points $u^{(0)} \in \ell_2(\Lambda,\RR^M)$,  $v^{(0)}$
                with $0\leq v_\lambda^{(0)} \leq \rho_\lambda$,\\
              & number $n_{\max}$ of outer iterations,\\
              & number of inner iterations $L_n$, $n=1,\hdots,n_{\max}$.\\
 Parameters: & $q \in [1,\infty]$, positive weights $(\theta_\lambda)$, $(\rho_\lambda)$, $(\omega_\lambda)$ with $\omega_\lambda \geq c > 0$,\\ 
& such that $\Phi^{(q)}$ and hence $J$
are convex, see Proposition \ref{q12oo}\\
 Output: & Approximation $(u^*,v^*)$ of the minimizer of $J^{(q)}_{\theta,\rho,\omega}$\\
    \end{tabular}
    \begin{tabbing}
      $u^{(0,0)} :=u^{(0)}$;\\
      {\em \texttt{for}} \= $n:=0$ {\em \texttt{to}} $n_{\max}$ {\em \texttt{do}}\\
      \>  {\em \texttt{for}} \= $m:=0$ {\em \texttt{to}} $L_n$ {\em \texttt{do}} \\
      \> \> $u^{(n,m+1)}:=U^{(q)}_{v^{(n)},\omega} \left (u^{(n,m)}+T^*(g - T u^{(n,m)}) \right) ;$\\
      \> {\em \texttt{endfor}} \\
        \> $u^{(n+1,0)}:= u^{(n,L_n)};$\\
      \> $v^{(n+1)}:=\left( \left \{ \begin{array}{ll} \rho_\lambda - \frac{1}{2 \theta_\lambda} \|{u^{(n+1,0)}}_\lambda\|_q, & \|{u^{(n+1,0)}}_\lambda\|_q < 2 \theta_\lambda \rho_\lambda\\
0,& \text{ otherwise }.
\end{array}
\right . \right)_{\lambda \in \Lambda};$\\
      
      {\em \texttt{endfor}}\\
$u^*:= u^{(n_{\max},L_{n_{\max}})};$\\
$v^* := v^{(n_{\max})}.$
    \end{tabbing}
  }
\end{alg}

Observe that each (inner) iteration of the above algorithm involves an
application of $T^*T$ and of the thresholding operator $U^{(q)}_{v,\omega}$. The latter
can be applied fast. So if there is also a fast algorithm for the computation of $T^* T$
then each iteration can be done fast.

Our analysis ensures the (weak) convergence of this scheme only if the inner loop computes
exactly the minimizer of $J(u,v^{(n)})$ for fixed $v^{(n)}$, i.e., if $L_n = \infty$. Of course,
this cannot be numerically realized, so we need to analyze what happens if the inner loop makes a small
error in computing this minimizer. In other words, how large do we have to choose $n_{\max}$ 
and $L_n, n=1,\hdots,n_{\max}$ in order to ensure that we have approximately 
computed the minimizer $u^*,v^*$ within a given error tolerance?



\subsection{Error analysis and strong convergence of {\bf JOINTSPARSE}}

First of all we want to establish the convergence rate of the inner loop, i.e.,
the iterative thresholding algorithm of the previous Section.

\begin{proposition} 
\label{err1}
Assume that $\omega_\lambda \geq \gamma > 0$ for all $\lambda \in \Lambda$
(implying that $K(u) = K^{(q)}_{v,\omega}(u)$ is strictly convex) and $\|T\| < 1$.
Set $\alpha := (1+\gamma)^{-1} \| I - T^* T\| < 1$.
Then the iterative thresholding algorithm
$$
 u^{(n,m+1)}\,:=\,U^{(q)}_{v^{(n)},\omega} \left (u^{(n,m)}+T^*(g - T u^{(n,m)}) \right),
$$
converges linearly 
\begin{equation}\label{err_est}
        \|u^{(n,\infty)} - u^{(n,m+1)}\|_2 \,\leq\, \alpha \|u^{(n,\infty)} - u^{(n,m)}\|_2.
\end{equation}
\end{proposition}
\begin{proof}
Note that 
$$
 u^{(n,\infty)}:=U^{(q)}_{v^{(n)},\omega} \left (u^{(n,\infty)}+T^*(g - T u^{(n,\infty)}) \right).  
$$
By non-expansiveness of $S_{v}^{(q)}$ (see Lemma \ref{nonexp} and its proof)
we obtain
\begin{align}
&\|u^{(n,\infty)} - u^{(n,m+1)}\|_2 \notag\\
= &\, \| U^{(q)}_{v^{(n)},\omega} \left (u^{(n,\infty)}+T^*(g - T u^{(n,\infty)}) \right) -U^{(q)}_{v^{(n)},\omega} \left (u^{(n,m)}+T^*(g - T u^{(n,m)}) \right)\|_2\notag\\
\,= & \left(\sum_{\lambda \in \Lambda} (1+\omega_\lambda)^{-2} \|S_{v_\lambda}^{(q)}((u^{(n,\infty)}+T^*(g - T u^{(n,\infty)})_\lambda) -  S_{v_\lambda}^{(q)}((u^{(n,m)}+T^*(g - T u^{(n,m)})_\lambda)\|_2^2\right)^{1/2}\notag\\
\leq\, & \sup_{\lambda \in \Lambda} (1+\omega_\lambda)^{-1} \| (I-T^*T)(u^{(n,\infty)} -u^{(n,m)})\|_2
\,\leq\, (1+\gamma)^{-1}\|I-T^*T\|\, \|u^{(n,\infty)} -u^{(n,m)}\|_2 \notag\\
=&\, \alpha \|u^{(n,\infty)} -u^{(n,m)}\|_2. \notag
\end{align}
This establishes the claim.
\end{proof}

\begin{rem} Clearly, the error estimation in (\ref{err_est}) holds also if one is only
interested in analyzing the iterative thresholding algorithm from the last section
(i.e. without doing the outer iteration).
Then it might also be interesting to consider the case that $\omega = 0$. According
to what we have proven in the previous section the algorithm
still converges provided the weight $v$ is bounded away from zero. However, 
then the error estimation (\ref{err1}) has a useful 
meaning only if $\alpha = \|I - T^*T\| < 1$.
So this applies if $T^*T$ is boundedly invertible. For a usual inverse problem,
however, we will have a non-invertible $T$ or at least one with unbounded inverse resulting in $\|I-T^*T\| = 1$.
So in this case we only know that the algorithm converges, but an error estimate does not
seem to be available.
\end{rem}

For simplicity we restrict the following error analysis to the most interesting cases $q\in \{1,2,\infty\}$.
We first need the following technical result.
\begin{lemma}
\label{lipproj}
For $q \in \{1,2,\infty\}$ the projection $P^q_v$ onto the ball $B^q(v) \subset \mathbb{R}^M$ is a Lipschitz function with respect to $v \in \mathbb{R}_+$. In particular, we have
\begin{equation}
\label{lipProj}
\|P^q_v(x) - P^q_w(x)|\ell_2^M\| \leq L |v-w| \quad \mbox{ for all } x \in \mathbb{R}^M,
\end{equation}
where $L=1$ for $q=2$ and $L=M^{1/2}$ for $q\in \{1,\infty\}$.
\end{lemma}
\begin{proof}
Let us start with $q=2$.
By distinguishing cases it is not difficult to show that
$$
\|P^2_v(x) - P^2_w(x)|\ell_2^M\| \leq |v-w|.
$$
For $q=\infty$ we have
$P^{\infty}_v(x) = (p^{\infty}_v(x_\ell))_{\ell = 1}^M$ where for $y \in \RR$
\[
p^{\infty}_{v}(y) \,=\,  \left\{\begin{array}{ll} y & \mbox{ if } |y| \leq v,\\
y- \sign(y)(|y|-v) & \mbox{ otherwise.}\end{array} \right.
\]
Since $p^{\infty}_{v}$ can be interpreted as a projection onto the $\ell_2$ ball in dimension 1, we obtain that
$$
        |p^{\infty}_{v}(y)-p^{\infty}_{w}(y)|\leq|v-w|, 
$$
and
$$
\|P^\infty_v(x) - P^\infty_w(x)|\ell_2^M\| = \left (\sum_{\ell=1}^M |p^{\infty}_{v}(x_\ell)-p^{\infty}_{w}(x_\ell)|^2\right)^{1/2} \leq M^{1/2} |v-w|.
$$
The case $q=1$ requires a bit more effort. 
By Lemma \ref{explS} (c) we have the following. 
Let $x_{i_k}$ denote the reordering of the entries of $x$ by magnitude
as in Lemma \ref{explS}. Let $n \in \{1,\hdots,M\}$ be the largest index satisfying
$$
|x_{i_n}| \,\geq\, \frac{1}{n-1}\left(\sum_{k=1}^{n-1} |x_{i_k}| - v\right).
$$
Then
\begin{align}
(P_v^{1}(x))_{i_j} \,&=\, x_{i_j} - 
\frac{\sign(x_{i_j})}{n}\left(\sum_{k=1}^n |x_{i_k}| - v\right), \quad j=1,\hdots,n, \notag\\
(P_v^{1}(x))_{i_j} \,&=\, 0,\quad j=n+1,\hdots,M.\notag
\end{align}
Observe first that for all $x \in \mathbb{R}^M$ there exists $\varepsilon_0>0$ such that for all $0 <\varepsilon < \varepsilon_0$ the same $n \in \{1,\hdots,M\}$ is the largest index satisfying
$$
|x_{i_n}| \,\geq\, \frac{1}{n-1}\left(\sum_{k=1}^{n-1} |x_{i_k}| - (v+\varepsilon)\right).
$$
For $0<\varepsilon < \varepsilon_0$, a simple computation yields 
\[
\frac{(P_{v+\varepsilon}(x)-P_v(x))_{i_j}}{\varepsilon} \,=\, \left\{
\begin{array}{ll} \frac{\sign(x_{i_j})}{n} & \mbox{ for } j=1,\hdots,n,\\
0 & j=n+1,\hdots,M.
\end{array}\right.
\]
This means that the map $v \rightarrow P^1_v(x)$ is right-differentiable, i.e., the limit 
$$
(P^1_v(x))'_+ = \lim_{\varepsilon \rightarrow 0_+} \frac{P^1_{v+\varepsilon}(x) - P^1_{v}(x)}{\varepsilon} 
$$
exists in $\mathbb{R}^M$. Moreover, it also follows that
\begin{equation}
\label{dxder}
\|(P^1_v(x))'_+\|_2 \,=\, \sqrt{n} \leq \sqrt{M}. 
\end{equation}
To conclude the proof we use the following standard result.
\begin{lemma}
Let $f:\mathbb{R} \rightarrow \mathbb{R}^M$ and $\varphi:\mathbb{R} \rightarrow \mathbb{R}$ be two continuous and right differentiable functions such that
$$
        \|f'_+(v)\| \leq \varphi'_+(v), 
$$
for all $v \in \mathbb{R}$. Then 
$$
        \| f(v) - f(w) \| \leq \varphi(v) - \varphi(w), \quad \mbox{ for all } v \geq w.
$$
\end{lemma}
According to the notation of this latter lemma, let us set $f(v) = P^1_v(x)$ and $\varphi(v)= M^{1/2} v$. Since $P^1_v(x)$ is a continuous function with respect to $v$ (in fact this is true for any projection onto convex sets, see \cite{rowe98}), by \eqref{dxder} and an application of the lemma we conclude that
$
        \| P^1_v(x) -  P^1_w(x) \| \leq M^{1/2} |v-w|.
$
\end{proof}

Observe that the strict convexity of $\Phi^{(q)}(u,v)$ is equivalent
to $\theta_\lambda \omega_\lambda > \kappa/4$, see Proposition \ref{q12oo}. 
In the following Proposition
we require the slightly stronger condition that $\theta_\lambda \omega_\lambda$
is bounded strictly away from $\kappa/4$, at least for $q=1,2$.

\begin{proposition}
\label{err2} Let $q \in \{1,2,\infty\}$. Assume that $\theta_\lambda \omega_\lambda \geq \sigma > \phi_q/4$  for all $\lambda \in \Lambda$, where $\phi_1 = M$, $\phi_2 = 1$, $\phi_\infty = \sqrt{M}$,
implying that $\Phi^{(q)}(u,v)$ and $J(u,v)$ are strictly convex, see Proposition \ref{q12oo}.
Moreover, let us assume that $\omega_\lambda \geq \gamma>0$ for all $\lambda \in \Lambda$.
Suppose $\|T\| < 1$ resulting in $\|I-T^*T\| \leq 1$. 
Set
\begin{equation}\label{def_beta}
\beta \,:=\, \sup_{\lambda \in \Lambda}\, \frac{\phi_q}{4\theta_\lambda\omega_\lambda + 4\theta_\lambda(1-\|I - T^*T\|)} \,\leq\, \frac{\phi_q}{4\sigma} \,<\, 1. 
\end{equation}
Then for each $n \in \mathbb{N}$ one has the following error estimate
\begin{equation*}
\|u^{(n,\infty)} - u^*|\ell_2(\Lambda,\RR^M)\| \,\leq\, \beta \| u^{(n,0)} - u^*|\ell_2(\Lambda,\RR^M)\|.
\end{equation*}
\end{proposition}
\begin{proof}
Let us consider the $n$-th iteration of the outer loop.
We have
\begin{equation*}
u^{(n,\infty)} \,=\, U_{v^{(n)},\omega}^{(q)}(\underbrace{u^{(n,\infty)} + T^*(g - T u^{(n,\infty)})}_{:=y^{(n)}}).
\end{equation*}
By the weak convergence of the double-minimization algorithm, also the minimum solution $u^*$ satisfies a similar relation,
\begin{equation*}
u^* \,=\, U_{v^*,\omega}^{(q)}(\underbrace{u^* + T^*(g - T u^*)}_{:=y^*}).
\end{equation*}
Recall that
$U_{v,\omega}^{(q)}(y)_\lambda := (1+\omega_\lambda)^{-1}S^{(q)}_{v_\lambda}(y_\lambda)$.
By non-expansiveness of $S^{(q)}_v$ (Lemma \ref{nonexp}) we have
\begin{eqnarray*}
\|u^{(n,\infty)}_\lambda - u^*_\lambda\|_2 &\leq & \| U_{v^{(n)},\omega}^{(q)}(y^{(n)})_\lambda- U_{v^{(n)},\omega}^{(q)}(y^*)_\lambda\|_2 + \| U_{v^{(n)},\omega}^{(q)}(y^*)_\lambda - U_{v^*,\omega}^{(q)}(y^*)_\lambda \|_2\\
&=& (1+\omega_\lambda)^{-1} \|y^{(n)}_\lambda - y^*_\lambda \|_2 + \| U_{v^{(n)},\omega}^{(q)}(y^*)_\lambda - U_{v^*,\omega}^{(q)}(y^*)_\lambda \|_2\\
&\leq & (1+\omega_\lambda)^{-1}\|I -T^*T\|\, \|u^{(n,\infty)}_\lambda - u^*_\lambda\|_2 + \| U_{v^{(n)},\omega}^{(q)}(y^*)_\lambda - U_{v^*,\omega}^{(q)}(y^*)_\lambda \|_2.
\end{eqnarray*}
This implies
\begin{align}
\label{firstestim}
\|u^{(n,\infty)}_\lambda - u^*_\lambda\|_2 \,&\leq\, \left(1-(1+\omega_\lambda)^{-1}\|I-T^*T\|\right)^{-1} \| U_{v^{(n)},\omega}^{(q)}(y^*)_\lambda - U_{v^*,\omega}^{(q)}(y^*)_\lambda \|_2\notag\\
& =\, \left(1+\omega_\lambda - \|I-T^*T\|\right)^{-1} 
\|S_{v_\lambda^{(n)}}(y_\lambda^*) - S_{v_\lambda^*}(y_\lambda^*)\|_2.\notag
\end{align}
Recall from Lemma \ref{lem_thresh1} that
$S^{(q)}_{v}{(2)}(x) \,=\, x - P_{v/2}^q(x)$ where $P^q_{v/2}$ denotes the orthogonal projection
of $x$ onto the $\ell_q$-ball of radius $v/2$. By Lemma \ref{lipproj} we have that for any
$z \in \RR^M$
\[
\|S^{(q)}_{v_\lambda^{(n)}}(z) - S^{(q)}_{v_\lambda^*}(z)\|_2 \,=\, 
\|P^q_{v_\lambda^*/2}(z) - P^q_{v_\lambda^{(n)}/2}(z)\|_2 \,\leq\, \frac{L}{2} |v_\lambda^{(n)} - v_\lambda^*|.
\]
So $S^{(q)}_v(z)$ is also Lipschitz in $v$. 
Let us recall that
$$
 v_\lambda^{(n)}:=\left \{ \begin{array}{ll} \rho_\lambda - \frac{1}{2 \theta_\lambda} \|u^{(n,0)}_\lambda\|_q, & \|u^{(n,0)}_\lambda\|_q < 2 \theta_\lambda \rho_\lambda\\
0,& \text{ otherwise }.
\end{array}
\right . 
$$
and
$$
 v_\lambda^*:=\left \{ \begin{array}{ll} \rho_\lambda - \frac{1}{2 \theta_\lambda} \|u^*_\lambda\|_q, & \|u^*_\lambda\|_q < 2 \theta_\lambda \rho_\lambda\\
0,& \text{ otherwise }.
\end{array}
\right . 
$$
By distinguishing cases we can show that
$$
        |  v_\lambda^{(n)} -  v_\lambda^* | \leq \frac{1}{2 \theta_\lambda} \left | \|u^{(n,0)}_\lambda\|_q - \|u^*_\lambda\|_q \right | \,\leq\,  \frac{1}{2 \theta_\lambda} \|u^{(n,0)}_\lambda - u^*_\lambda\|_q \leq  \frac{R}{2 \theta_\lambda} \|u^{(n,0)}_\lambda - u^*_\lambda\|_2,
$$      
where $R=1$ for $q \in \{2,\infty\}$ and $R=M^{-1/2}$ for $q=1$.
Pasting the pieces together yields 
\[
\|u_\lambda^{(n,\infty)} - u^*_\lambda\|_2
\,\leq\, \frac{\phi_q}{4\theta_\lambda\left(\omega_\lambda + 1 -\|I-T^*T\|\right)} \|u^{(n,0)}_\lambda - u^*_\lambda\|_2
\,\leq\, \beta \|u^{(n,0)}_\lambda - u^*_\lambda\|_2.
\]
Summation over $\lambda \in \Lambda$ completes the proof.
\end{proof}

Let us combine the previous two results to obtain the error estimation for the finite algorithm, i.e., for $L_n<\infty$.
\begin{theorem} Make the same assumptions as in Propositions \ref{err1} and \ref{err2}.
Choose $L_{n}$ such that
\[
\delta_n \,:=\,  \left(\alpha^{L_{n}}(1+\beta) + \beta\right) \leq \delta < 1\quad\mbox{ for all } n \in \bN.
\]
(This is possible since $\alpha,\beta < 1$). Then we have linear convergence of our algorithm, i.e.,
\begin{equation*}
\|u^{(n,0)} - u^*\|_2 \leq \delta_n \|u^{(n-1,0)} - u^*\|_2.
\end{equation*}
\end{theorem}
\begin{proof}
Using Proposition \ref{err1} and Proposition \ref{err2} we get 
\begin{eqnarray*} 
\|u^{(n,0)} - u^*\|_2 &\leq& \| u^{(n,0)} - u^{(n-1,\infty)}\|_2 + \| u^{(n-1,\infty)} - u^*\|_2\\
&\leq& \alpha^{L_{n-1}} \|u^{(n-1,0)} - u^{(n-1,\infty)}\|_2 + \beta \| u^{(n-1,0)} - u^*\|_2 \\
& \leq & \alpha^{L_{n-1}} \left(\|u^{(n-1,0)} - u^*\|_2 + \|u^* - u^{(n-1,\infty)}\|_2\right)
+ \beta\| u^{(n-1,0)} - u^*\|_2 \\
&\leq & \alpha^{L_{n-1}}\left(\|u^{(n-1,0)} - u^*\|_2 + \beta\|u^{(n-1,0)} - u^*\|_2\right) 
+ \beta \| u^{(n-1,0)} - u^*\|_2\\
&\leq & \left(\alpha^{L_{n-1}}(1 + \beta) + \beta\right) \|u^{(n-1,0)} - u^*\|_2.
\end{eqnarray*}
This concludes the proof.
\end{proof}

\begin{rem}
The last theorem shows that it is possible to choose the number $L_n$ of inner iterations
constant with respect to $n$. 
\end{rem}

\subsection{Strong convergence of the double-minimization algorithm}

Finally, we can establish the strong convergence of the double-minimization algorithm and conclude the full proof of  Theorem \ref{thm_conv}.

\begin{corollary} Under the assumptions of Proposition \ref{err2}, if the minimizer 
of $J(u,v^{(n)})$ for fixed $v^{(n)}$ 
could be computed exactly, i.e., $L_n=\infty$ for all $n \in \mathbb{N}$, 
then the outer loop converges with exponential rate, and we have
\begin{equation*}
\|u^{(n)} - u^*\|_2 \,\leq\,\beta^n \| u^{(0)} - u^*\|_2,
\end{equation*}
where we have denoted here $u^{(n)}:=u^{(n,\infty)}$. Moreover, the sequence $v^{(n)}$ converges
componentwise and $v^{(n)}-v^*$ converges to $0$ strongly in $\ell_{2,\theta}(\Lambda)$.
\end{corollary}
\begin{proof}
The first part of the statement is a direct application of Proposition \ref{err2}. It remains to 
show that $v^{(n)}-v^*$ converges to $0$ strongly in $\ell_{2,\theta}(\Lambda)$.
Using that all norms on $\RR^M$ are equivalent it follows that 
\begin{align}
4 \sum_{\lambda \in \Lambda} \theta_\lambda^2|v^{(n)} - v^{*}|^2 \,&=\, 
\sum_{\lambda \in \Lambda} |\|u^{*}_\lambda\|_q - \|u^{(n)}_\lambda\|_q|^2
\,\leq\, \sum_{\lambda \in \Lambda} \|u^{*}_\lambda - u^{(n)}_\lambda\|_q^2\notag\\
&\leq\, C \sum_{\lambda \in \Lambda} \|u^{*}_\lambda - u^{(n)}_\lambda\|_2^2
\,=\, C \|u^{*} - u^{(n)}|\ell_{2}(\Lambda,\RR^M)\|^2. \notag
\end{align}
Thus, $v^{(n)}-v^{(\infty)}$ converges also strongly in $\ell_{2,\theta}(\Lambda)$.
\end{proof}

\section{Color image reconstruction}

With this section we illustrate the application of the algorithms for color image recovery. The scope 
is to furnish a qualitative description of the behavior of the scheme. In a subsequent work we plan to 
provide a finer quantitative analysis 
in the context of {\it distributed compressed sensing} \cite{babadusawa05}.

We begin by illustrating an interesting real-world problem occurring in art restoration. 
On $11^{th}$ March 1944, a group of bombs launched from an Allied airplane hit the famous Italian Eremitani's 
Church in Padua, destroying it together with the inestimable frescoes by Andrea Mantegna {\it et al.} contained 
in the Ovetari Chapel. Details on ``the state of the art'' can be found in \cite{FT1,FT2}. In 1920 a collection 
of high quality gray level pictures of these frescoes has been made by Alinari. The only color images of the 
frescoes are dated to 1940, but unfortunately their quality (i.e., the intrinsic resolution of the printouts) 
is much lower, see Figure \ref{resultfresco}.
Inspired by the fresco application, we model the problem of the recovery of a high resolution 
color image from a low resolution color datum and a high resolution gray datum. We will implement 
the solution to the model problem as a non-trivial application of the algorithms we have presented in this paper.

\subsection{Color images, curvelets, and joint sparsity}

Let us assume that the color images are encoded into YIQ channels.
The Y component 
represents the luminance information (gray level), while I and Q give
the chrominance information. 
Of course, one may also choose a different encoding system, e.g., RGB or CMYK. 
Clearly, the color image $f=(f_{1},f_{2},f_{3})$ can be represented as a 3-channel signal.  
In order to apply our algorithm, we need to fix a frame for which we 
can assume color images being jointly sparse.

It is well-known that {\it curvelets} \cite{cado04} are well-suited for sparse approximations 
of curved singularities. 
A natural image can in 
fact be modelled as a function which is piecewise smooth except on a discontinuity set, the latter 
being described as the union of rectifiable curves.  
Moreover, there are fast algorithms available for the 
computation of curvelet coefficients of digital images \cite{cadedoyi50}.

In the following, let us assume that a color image $f$ is encoded into a vector of curvelet 
coefficients $(u^\ell_\lambda)_{\lambda \in \Lambda}^{\ell=1,2,3}$. 
The image can be reconstructed by the synthesis formula
$$
        f \, = \, (F u^\ell)_{\ell=1,2,3}:= 
\left (\sum_{\lambda \in \Lambda} u^\ell_\lambda \psi_\lambda \right)_{\ell=1,2,3},
$$
where $\{\psi_\lambda: \lambda \in \Lambda\}$ is the collection of curvelets.
The index $\lambda$ consists of 3 different parameters, $\lambda=(j ,p, k)$, where
$j$ corresponds to scale, $p$ to a rotation, and $k$ 
to the spatial location of the curvelet $\psi_\lambda$. 
We do not enter in further details, especially of the discrete and numerical implementation, 
which one can find in \cite{cado04,cadedoyi50}.

Let us instead observe that significant curvelet coefficients $u_\lambda=(u^1_\lambda, u^2_\lambda, u^3_\lambda)$ 
will appear simultaneously at the same $\lambda \in \Lambda$ for all the channels, as soon as the corresponding 
curvelet overlaps with a (curved) singularity (appearing simultaneously in all the channels), and is approximately 
tangent to it. This justifies the joint sparsity 
assumption for color images with respect to curvelets.

 \begin{figure}[h]
  \centering 
\includegraphics[width=.45\textwidth]{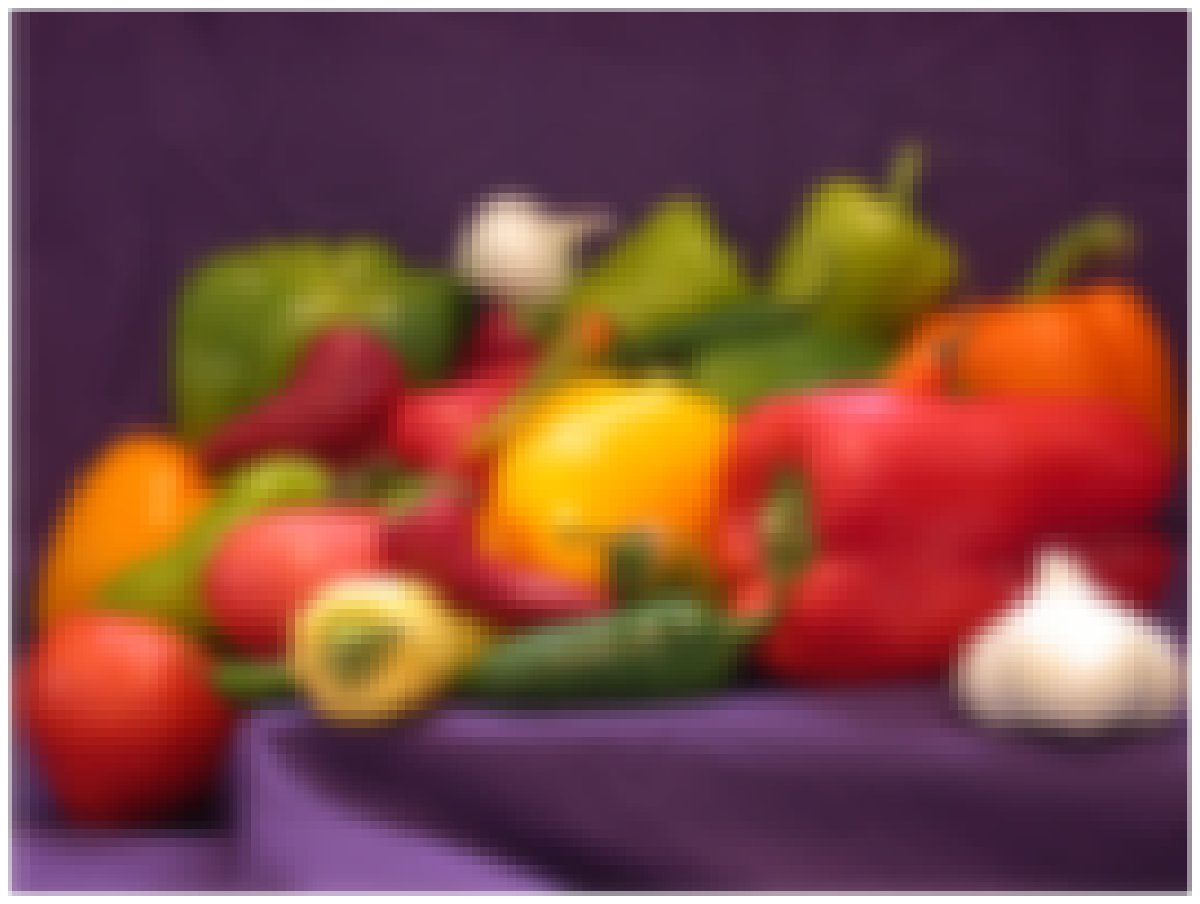}
 \includegraphics[width=.45\textwidth]{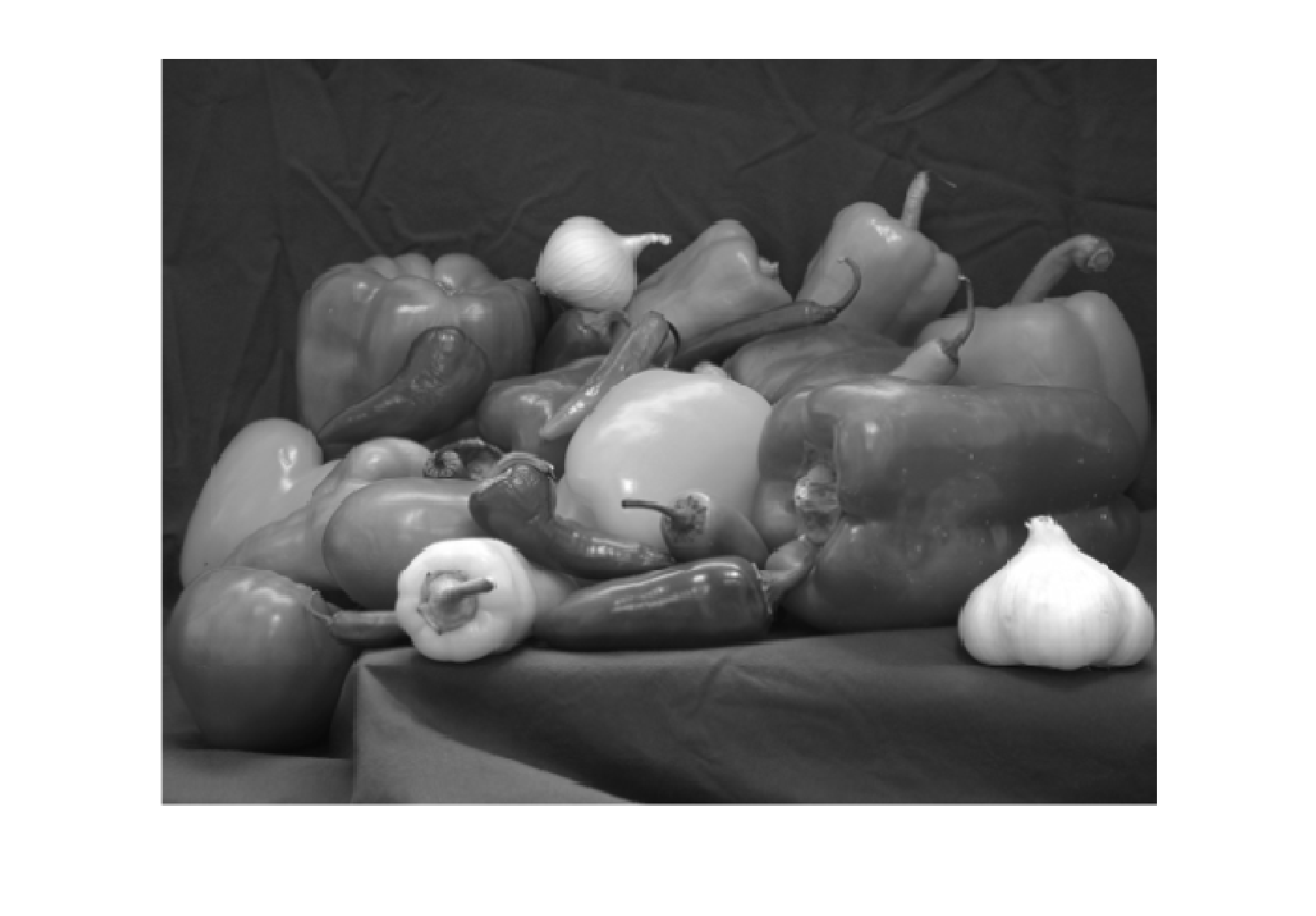}
  \caption{{\em Left:} The low-resolution color image is here presented after Gaussian filtering and downsampling. In the numerical experiments the I and Q channels are used. {\em Right:} The high resolution gray level image encodes several morphological information useful to recover the high resolution color image.}
  \label{modelproblem}
\end{figure}
\subsection{The model of the problem}

The datum of our problem is a three-channel signal 
$g=(g_1,g_2,g_3) \in \ell_2(\mathbb{Z}_{N_0}^2,\mathbb{R}^3)$ 
where $g_i$, $i=2,3$, are the low resolution chrominance channels I and Q, 
and $g_1$ is the high resolution gray 
channel Y. 
We assume that $g$ was produced by
$g = T u$
where $u = (u^1,u^2,u^3)$ are the curvelet coefficients of the three channels
of the high resolution color image that we want to reconstruct. The operator $T=(T_{\ell,j})_{\ell,j=1,2,3}$ 
can be expressed by the matrix
\begin{equation}
T=\left ( \begin{array}{ccc}
F &0&0\\
0& A F &0\\
0&0& A F
\end{array}
\right).
\end{equation} 
Here, $A$ is the linear operator that transforms the high-resolution image into the low resolution
image.
In particular, $A$ can be taken as a convolution operator (with a Gaussian for instance)
followed by downsampling.
Eventually, we may assume a suitable scaling in order to make $\|T\| < 1$, and a different weighting
of the gray channel and the I,Q channel in the discrepency term. 
Since $A$ is not invertible, also the operator $T$ is not invertible, and the minimization of 
$\mathcal{T}(u)$ requires a regularization. Clearly, for this task we use the functional
$K$ defined in (\ref{def_Psi0}) or $J$ defined in (\ref{def_func_J}). 

\begin{figure}[h]
  \centering
\includegraphics[width=1\textwidth]{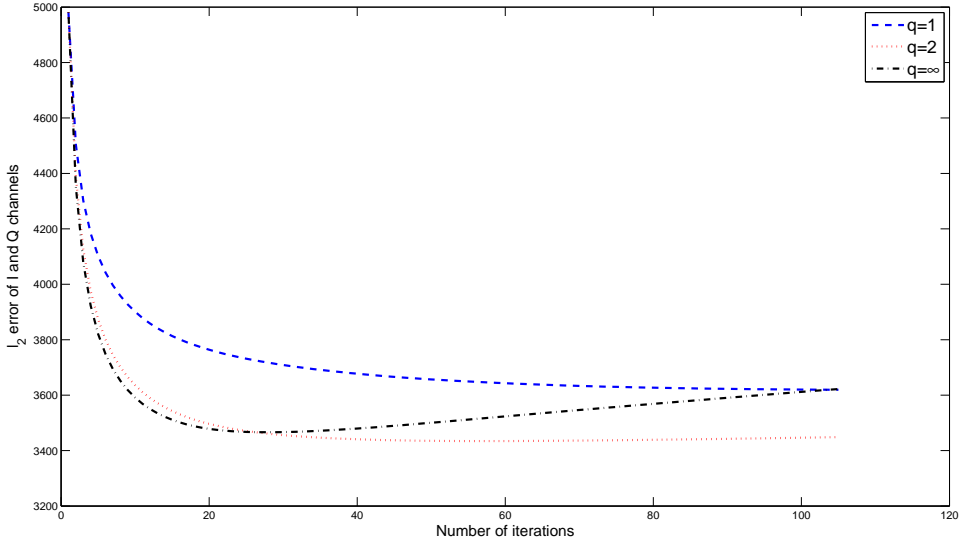}
  \caption{The $\ell_2$ error between the {\it original color image} and the iterations of the algorithm is shown for different values of $q=1,2,\infty$. We have considered $L_n=105$ and $n_{\max}=1$, the numbers of inner and outer iterations respectively. We have fixed here $\omega_\lambda=0$, $\theta_\lambda=10$, and $\rho_\lambda=20 \times 2^{-j}$.}
  \label{errorcomp1}
\end{figure}
 
\begin{figure}[h]
  \centering
\includegraphics[width=1\textwidth]{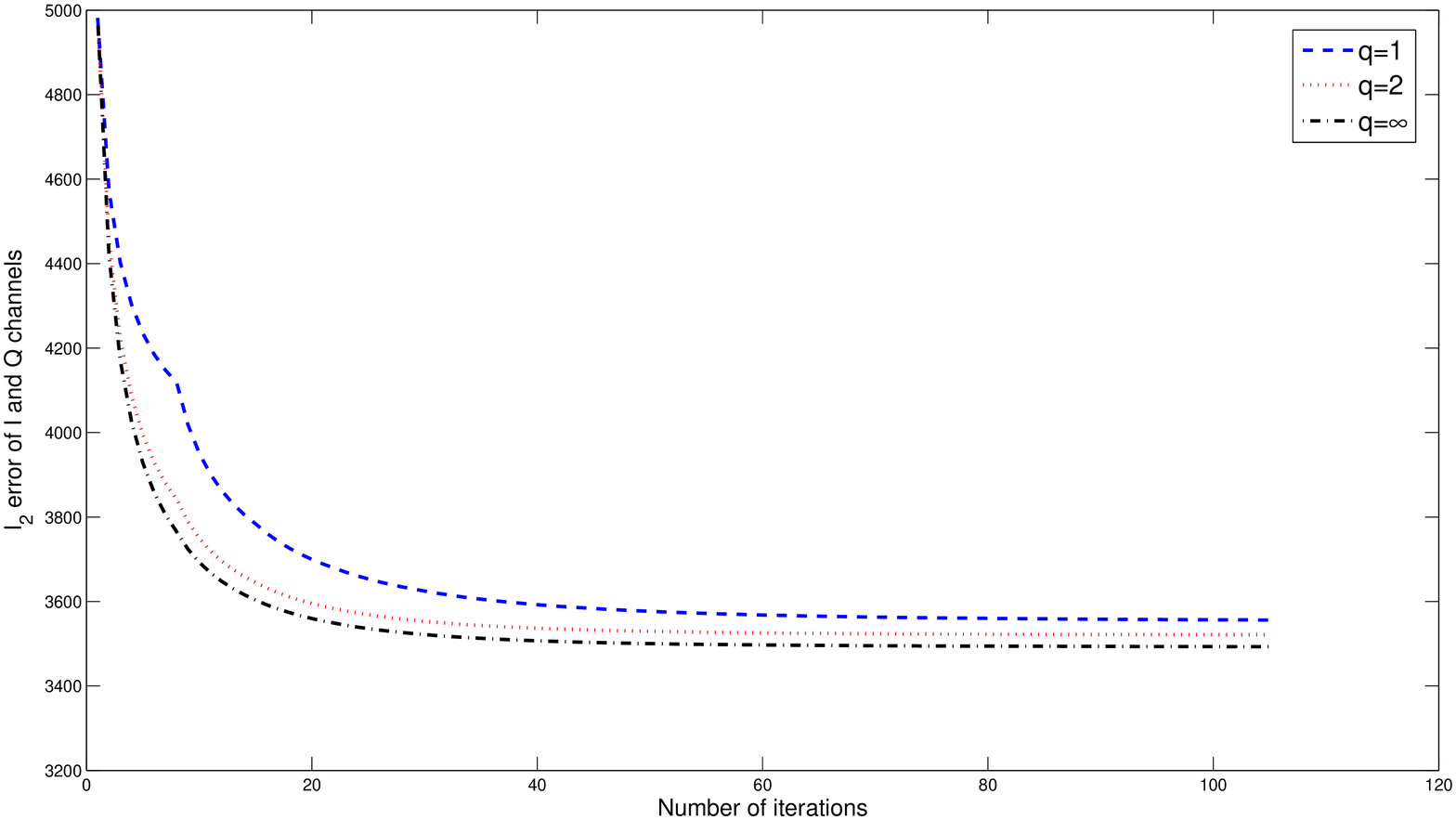}
  \caption{The $\ell_2$ error between the {\it original color image} and the iterations of the algorithm is shown for different values of $q=1,2,\infty$. We have considered $L_n=7$ and $n_{\max}=15$, the numbers of inner and outer iterations respectively. We have fixed here $\omega_\lambda=1/20$, $\theta_\lambda=10$, and $\rho_\lambda=20 \times 2^{-j}$.}
  \label{errorcomp2}
\end{figure}
\begin{figure}[h]
  \centering
\includegraphics[width=1\textwidth]{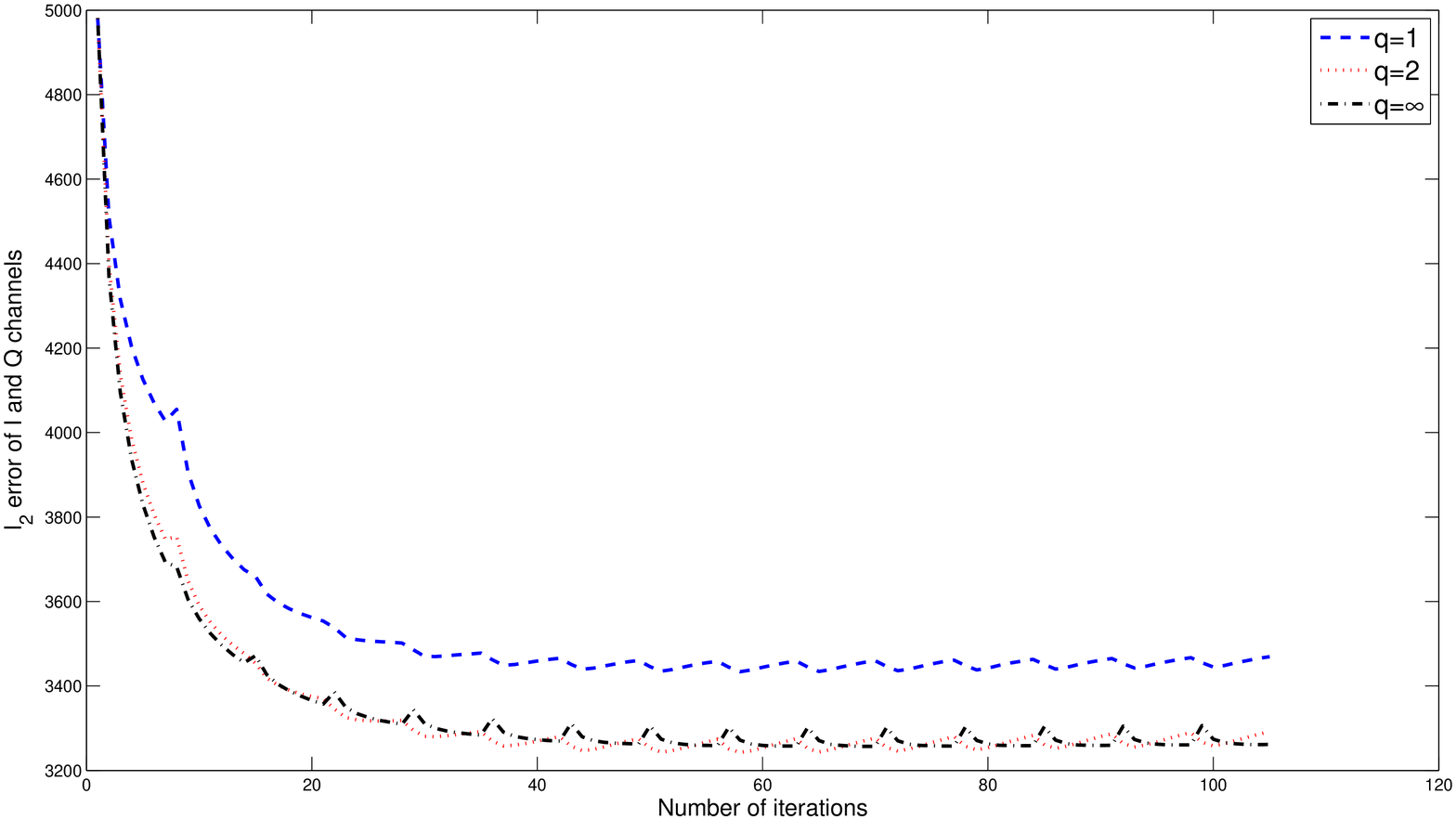}
  \caption{The $\ell_2$ error between the {\it original color image} and the iterations of the algorithm is shown for different values of $q=1,2,\infty$. We have considered $L_n=7$ and $n_{\max}=15$, the numbers of inner and outer iterations respectively. We have fixed here $\omega_\lambda=1/20$, $\theta_\lambda=10$, and $\rho_\lambda=20 \times 2^{-j}$. Further iterations of TV minimization are added in the outer loop to enforce edge enhancing.}
  \label{errorcomp3}
\end{figure}

\subsection{On the choice of the parameters}
\label{param}

What remains to clarify is the choice of the parameters $\omega_\lambda$, $\theta_\lambda$, 
and $\rho_\lambda$. The parameter $\omega_\lambda \geq \gamma >0$ has been introduced 
for the sole purpose to make $J$ strictly convex. 
A large value of this parameter actually produces an image $u$ which is significantly blurred and no information 
about edges is recovered. Thus, we rather put $\omega_\lambda = \gamma = \varepsilon >0$ small. 
Due to the convexity requirements (see Subsection 2.3), we select $\theta_\lambda \sim \frac{M}{\varepsilon}$.
The choice of $\rho_\lambda$ requires a deeper understanding of the information 
encoded by the curvelet coefficients.

Indeed, in \cite{cado04} it was observed that those curvelets 
that overlap with a discontinuity decay like $\|u_\lambda\|_q \lesssim 2^{-3/4 j}$ while
the others satisfy $\|u_\lambda\|_q \lesssim 2^{-3/2 j}$ (where $j$ denotes the scale).
Since we want to recover joint discontinuities
we may choose $\rho_\lambda:=\rho_{j,p,k}  \sim 2^{- j s}$ with $s \in [3/4, 3/2]$. By this choice and
by (\ref{min_v})
the locations $\lambda$ for which $v_\lambda =0$ will indicate a potential joint discontinuity.

Of course, this is just one possible choice of the parameters and further information 
might be extracted from the joint sparsity pattern indicated by $v$, by the use of different 
parameters. We believe that a deeper study of the characterization of the morphological 
properties of signals encoded by frames (e.g., curvelets and wavelets) is fundamental for the 
right choice of these parameters. We refer to \cite{ja04,ja05} for deeper insights in this 
direction,  concerning  fine properties of functions encoded by the distribution of 
wavelet coefficients.

\subsection{Numerical experiments}

According to the previous subsections, we illustrate here the application of {\bf JOINTSPARSE} 
for the recovery of a high resolution color image from a low resolution color datum 
and a high resolution gray datum. In Figure \ref{modelproblem} we illustrate the data of the problem. 
In this case the resolution of the color image has been reduced by a factor of 4 in each direction
by using a Gaussian filter 
and a downsampling. We have conducted several experiments for different choices 
of $q \in \{1,2,\infty\}$, with fixed parameters as indicated in Subsection \ref{param}. 
We have chosen $L_n=105$ and $n_{\max}=1$, as well as $L_n=7$ and $n_{\max}=15$ 
(the numbers of inner and outer iterations, respectively).  In the first case, only the 
minimization of $J(u,v^{(0)})$ with respect to $u$ has been performed, i.e., no iterative adaptation 
of the joint sparsity pattern indicated by $v$ occurred. In order to estimate the 
different behavior depending of the parameters above, we have evaluated at each 
iteration the $\ell_2$-error between the reconstructed I and Q color channels 
and the {\it original} I and Q color channels.  
Figures \ref{errorcomp1} and \ref{errorcomp2} indicate that the error decreases for increasing values of $q$. 
This means that the increased coupling due to the $q$-parameter is significant 
in order to improve the recovery. Recall that the choice $q=1$ does not induce 
any coupling between channels.

\begin{figure}[h]
\[\begin{array}{ccc}
  \centering
  
 \includegraphics[width=4cm]{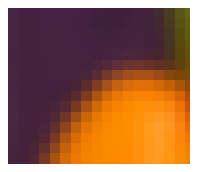} &  
  \includegraphics[width=4cm]{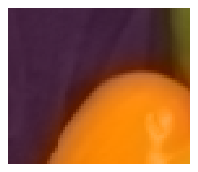} &  
  \includegraphics[width=4cm]{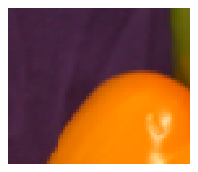}\\
\includegraphics[width=4cm]{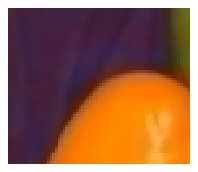} &  
  \includegraphics[width=4cm]{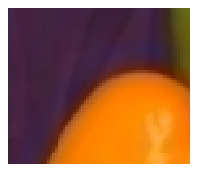} &  
  \includegraphics[width=4cm]{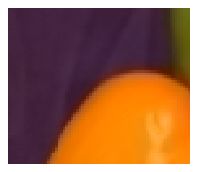}
\end{array}\]
\caption{
{\bf First row}. 
{\em Left:} Portion of the low resolution color image. 
    {\em Center:} Portion of the reconstructed color image by Gaussian interpolation and substitution of the $Y$ channel with the gray level datum. Evident color artifacts appear at edges.  {\em Right:} Portion of the original color image.
{\bf Second row}. 
{\em Left:} Portion of the reconstructed color image for $q=1$, $L_n=105$, and $n_{\max}=1$. 
    {\em Center:} Portion of the reconstructed color image for $q=\infty$, $L_n=7$, and $n_{\max}=15$.  {\em Right:} Portion of the reconstructed color image for $q=\infty$, $L_n=7$, $n_{\max}=15$, and TV minimization.}
\label{fig:compar}
\end{figure}

This coupling effect due to $q>1$ is even more evident in Figure \ref{errorcomp2}, 
where the adaptation of the weight $v$ occurs. The left and the central pictures in 
the second row of Figure \ref{fig:compar} show a reduced color distortion at edges, passing from the case $q=1$ 
(without coupling) to the case $q=\infty$ respectively, and consequently a better edge resolution. 
Nevertheless, the differences are not so remarkable. This is due to the fact that, although the 
functional $J$ promotes coupling at edges, it does not necessarily enforce a significant edge 
enhancement. Thus, we may modify the functional by adding an additional total variation 
constraint on the I and Q channels:
$$
        J_{\mbox{\small TV}}(u,v):= J(u,v) + \left(|F u^2|_{\mbox{\small TV}}+ |F u^3|_{\mbox{\small TV}}\right).
$$
The effect of this modification is to promote edge enhancing together with their simultaneous coupling through different channels. 
For the minimization of $J_{\mbox{\small TV}}$ we use a heuristic scheme as in \cite{elstqudo05}, by alternating iterations for the minimization of $J$ and for the minimization of $\left(|F u^2|_{\mbox{\small TV}}+ |F u^3|_{\mbox{\small TV}}\right)$, compare also \cite{dateve06}.
The corresponding results are shown in Figure \ref{errorcomp3} where the effect of the coupling (for the cases $q=2,\infty$) is further enhanced.
The right  picture in the second row of Figure \ref{fig:compar} shows the result of the reconstruction in this latter case. The edges are perfectly recovered.

\begin{figure}[h]
\[\begin{array}{ccc}
  \centering  
  \includegraphics[width=.345\textwidth]{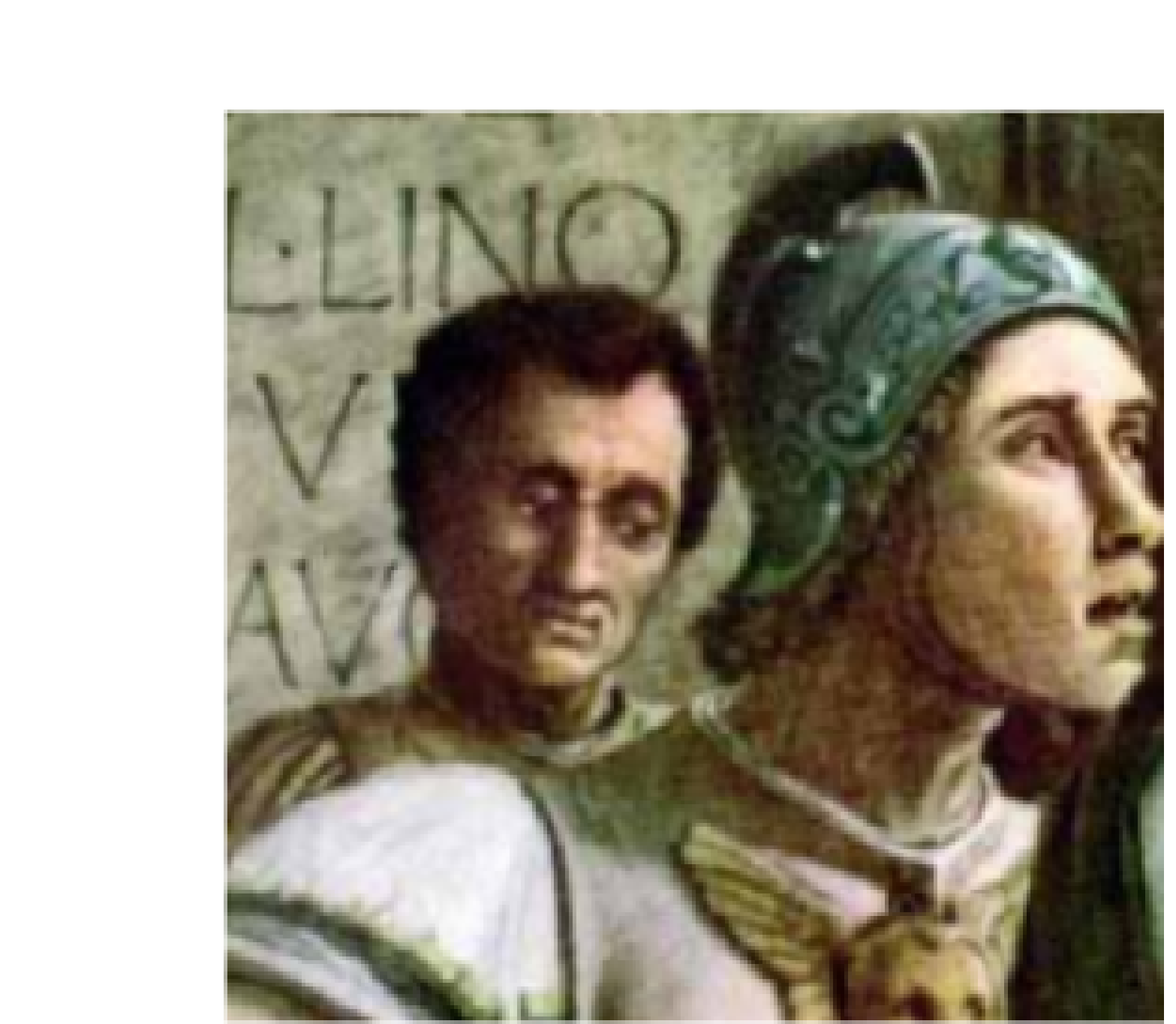}  &
  \includegraphics[width=.28\textwidth]{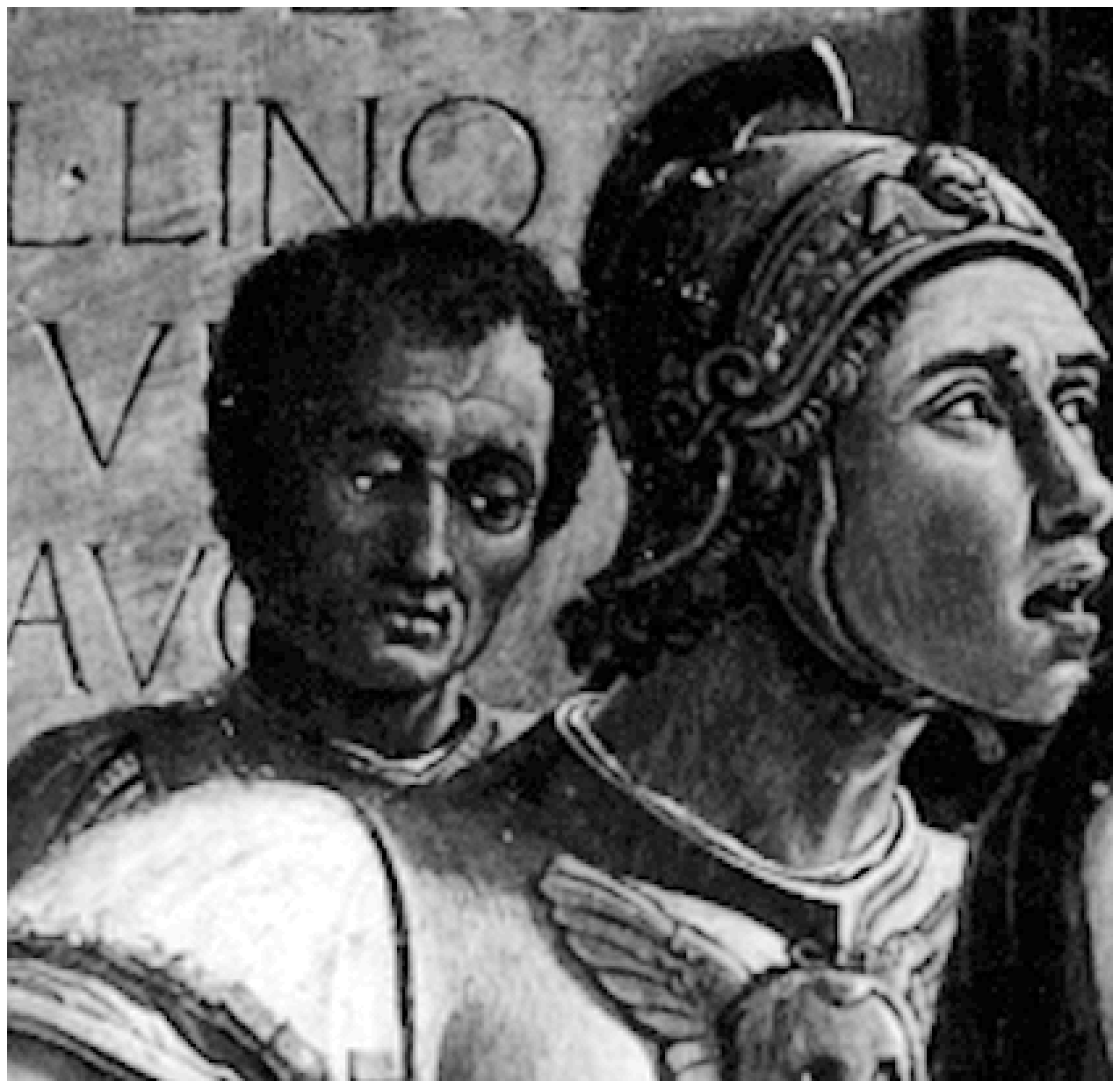}&
\includegraphics[width=.28\textwidth]{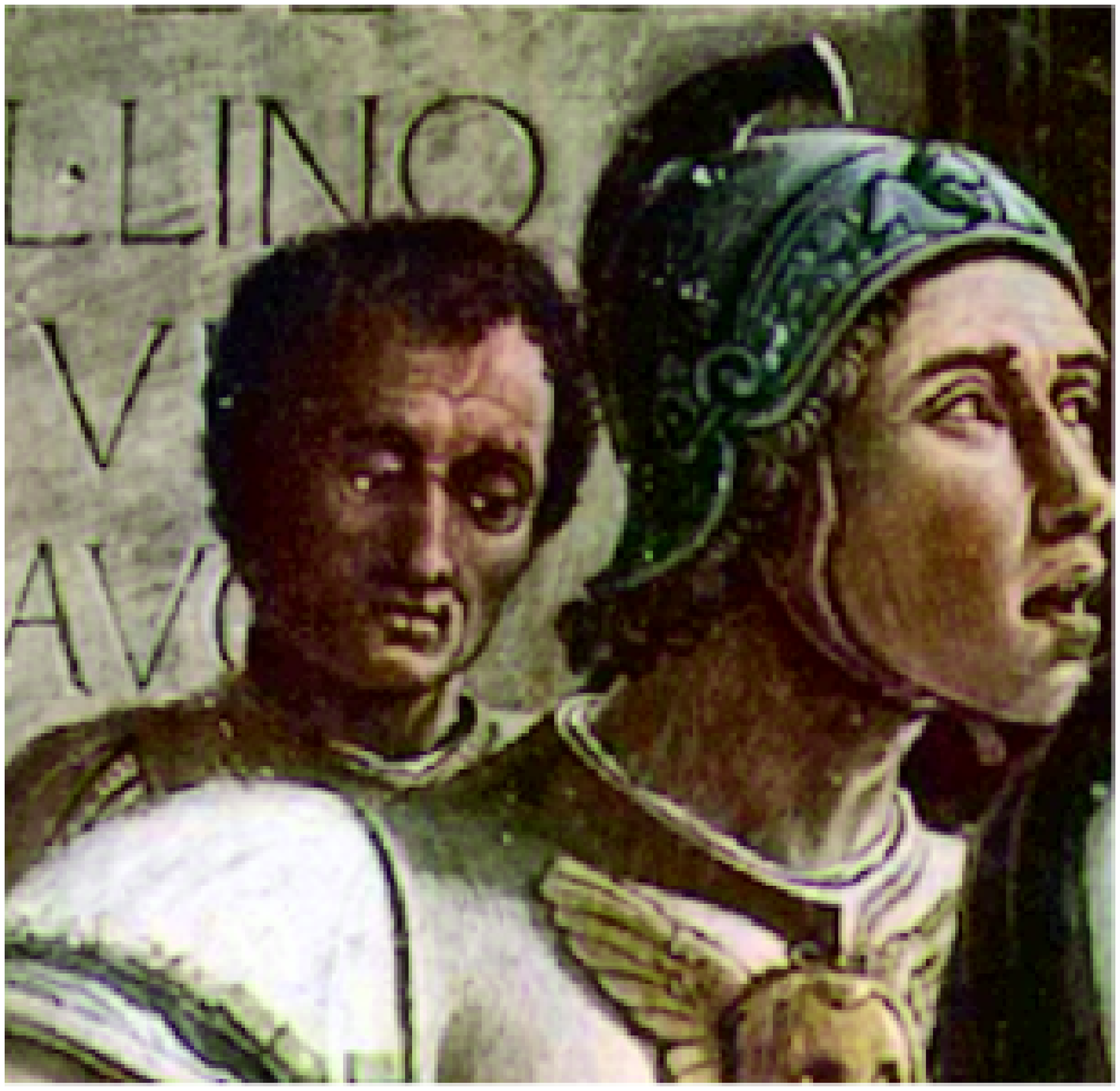}  
\end{array}\]
  \caption{{\em Left:}  Low quality color image of the fresco dated to 1940.
    {\em Center:} High quality gray image of the fresco dated to 1920. Some details are not visible in the color version. {\em Right:} The  reconstructed image after 6 outer iterations with 7 inner iterations each, for $q=\infty$. The final $Y$ channel is substituted with the high resolution gray level datum. The discountinuities are enhanced and no artifact colors appear.}
  \label{resultfresco}
\end{figure}

These numerical experiments confirm that the use of the joint sparsity measure 
$\Phi^{(q)}$ associated to the curvelet representation can improve significantly the quality 
of the reconstructed color image. Better results are achieved by choosing $q=\infty$ and 
by the adaptive choice of weights as indicators of the sparsity pattern. 
Further improvements can be achieved by channelwise edge enhancing, e.g.,  via total variation minimization.
An application to the real case of the art frescoes is illustrated in Figure \ref{resultfresco}.

\section{Final Remarks}


1. If the index set $\Lambda$ is infinite then $T^* T$ is represented as a 
biinfinite matrix and thus its evaluation might not be exactly numerically implementable. In a subsequent work 
we will consider the case $\# \Lambda = \infty$ and the treatment of sparse (approximate) 
evaluations of biinfinite matrices in order to realize fast and convergent schemes also in this
situation, compare also \cite{S,DFR,DFRSW}.

2. To exploit the optimal performance of the scheme, an extensive campaign of numerical experiments 
should be conducted in order to further refine the choice of parameters. 
It is also crucial to investigate the deeper relations among the parameter $\rho_\lambda$, 
the {\it multifractal analysis} as, e.g., in \cite{ja04}, and morphological image analysis. 
In particular, the parallel between the functional $J$ and the $\Gamma$-approximation of the 
Mumford-Shah functional by Ambrosio and Tortorelli \cite{amto90,brkiso03} is suggestive:
$$
        F_\varepsilon(u,v) := \underbrace{\int_\Omega (u - g)^2 dx}_{\sim \mathcal{T}(u)} + \underbrace{\int_\Omega v^2 f(\nabla u) dx}_{\sim \sum_\lambda v_\lambda \|u\|_q} + \underbrace{\int_\Omega \varepsilon |\nabla u|^2 dx}_{\sim \sum_\lambda \omega_\lambda \|u\|_2^2}  + \underbrace{\int_\Omega \left (\varepsilon |\nabla v|^2 + \frac{1}{4 \varepsilon} (1-v)^2\right ) dx}_{\sim \sum_\lambda \theta_\lambda(\rho_\lambda - v_\lambda)^2},
$$
where $f$ is a suitable polyconvex function, e.g.,  
$f(\nabla u) = (|\nabla u|^2 + |u_x \times u_y|) = (|\nabla u|^2 +| \mbox{adj}_2(\nabla u)|)$,  $\mbox{adj}_2(A)$ 
is the matrix of all $2 \times 2$ minors of $A$. The minimization of this term enforces that derivatives 
of different channels are large only in the same directions.
According to the specific choices of $\rho_\lambda$ to indicate the discontinuity set of $u$, and 
for  $\omega_\lambda = \varepsilon$ and $\theta_\lambda = \frac{1}{4 \varepsilon}$, we may investigate 
the behavior of the functional $J$ for $\varepsilon \rightarrow 0$ and its relation with the Mumford-Shah 
functional. The term $\sum_\lambda \frac{1}{4 \varepsilon} (\rho_\lambda - v_\lambda)^2$ essentially counts 
the number of curvelets that, from a certain scale $j$ on such that $2^{-j} \sim \varepsilon$, 
do overlap with the discontinuity 
set and are nearly tangent to the singularity. 
We conjecture that for $\varepsilon \rightarrow 0$ and for 
a rectifiable curved discontinuity, this term estimates the length of the discontinuity.

3. While we were finishing this paper, we have been informed by G. Teschke of the 
results in \cite{dateve06}. In this manuscript the authors consider linear inverse problems 
where the solution is assumed to fulfill some general 1-homogeneous convex constraint. 
They develop an algorithm that amounts to a projected Landweber iteration and that provides 
an iterative approach to the solution of this inverse problem. In particular for the case 
$\omega=(\omega_\lambda)_\lambda = 0$, some of our results stated in Section 4 can be 
reformulated in this more general setting and therefore derived from \cite{dateve06}. 
However, for $\omega \neq 0$ the sparsity measure $\Psi^{(q)}_{v,\omega}$ as in \eqref{def_Psi} 
is not 1-homogeneous and the elaborations in Section 4 are needed. Moreover, for the 
relevant cases $q=1,2,\infty$, we express explicitly the projection $P^{q'}_{v/2}$. 
Due to their generality, the results in \cite{dateve06} do not provide 
concrete recipes to compute such projections.

\section{Conclusion}

We have investigated joint sparsity measures with respect to frame 
expansions of vector valued functions. These sparsity measures generalize 
approaches valid for scalar functions and take into account common sparsity 
patterns through different channels. We have analyzed linear inverse problems 
with joint sparsity regularization as well as their efficient numerical 
solution by means of a novel algorithm based on thresholded 
Landweber iterations.
We have provided the convergence analysis for a wide range of parameters.
The role of the joint sparsity measure is twofold: to tighten the 
characterization of solutions of interest and to extract significant 
morphological properties which are a common feature of all the channels. 
By numerical applications in color image restoration, we have shown that 
joint sparsity significantly outperforms uncoupled constraints. 
We have presented the results of an application to a relevant 
real-world problem in art restoration.
The wide range of applicability of our approach includes 
several other problems with coupled vector valued solutions, 
e.g., neuroimaging and distributed compressed sensing.

\bibliographystyle{siam.bst}
\bibliography{FornasierRauhutBib}



\end{document}